\documentstyle[12pt]{article}
\begin{document}
\parindent 16pt
\title{\Large\bf A pre-order principle and  set-valued Ekeland variational  principle$^1$}
\setcounter{footnote}{1}
 \footnotetext{ This work was supported by
the National Natural Science Foundation of China (10871141).

 E-mail address: qjhsd@sina.com,  jhqiu@suda.edu.cn}

\author{  {Jing-Hui Qiu
}\\
{\footnotesize\sl School of Mathematical Sciences,  Soochow
University, Suzhou 215006, P. R. China} \\
}
\date{}
\maketitle
\begin{center}
\begin{minipage}{124mm}
\vskip 0.5cm {{\bf Abstract.}  \  We establish a pre-order
principle. From the principle, we obtain a very general set-valued
Ekeland variational principle, where the objective function is a
set-valued map taking values in a quasi ordered linear space and the
perturbation contains a family of set-valued maps satisfying certain
property.  From this general set-valued Ekeland variational
principle, we deduce a number of particular versions of  set-valued
Ekeland variational principle,
which include many known Ekeland variational principles, their  improvements   and some new results.\\

}\\


\end{minipage}
\end{center}
\vskip 1cm \baselineskip 18pt

\section*{ \large\bf 1. Introduction  }

\hspace*{\parindent}  In 1972, Ekeland [13] (see also [14, 15]) gave
a variational principle, now known as Ekeland variational principle
(for short, EVP), which says that for any lower semicontinuous
function $f$ bounded from below on a complete metric space, a
slightly perturbed function has a strict minimum.  In the last four
decades, the famous EVP emerged as one of the most important results
of nonlinear analysis and it has significant applications in
optimization, optimal control theory, game theory, fixed point
theory, nonlinear equations, dynamical systems, etc; see for example
[3, 10, 14, 15, 21, 36, 51]. Motivated by this wide usefulness, many
authors have been interested in extending EVP to vector-valued maps
or set-valued maps with values in a vector space quasi ordered by a
convex cone, see, for example,  [2, 4-10, 12, 16, 17, 21-25, 27-31,
33-35, 40-42, 44, 45,  47, 48] and the references therein.

Recently, there are many new and interesting results of EVP for
set-valued maps. Here, we only mention some results which are
related to this paper. In [24], Ha introduced a strict minimizer of
a set-valued map by virtue of Kuroiwa's set optimization criterion
(see [32]). Using the concept of cone extensions and
Dancs-Hegedus-Medvegyev theorem (see [11]) she established a new
version (see [24, Theorem 3.1]) of EVP for set-valued maps, which is
expressed by the existence of a strict minimizer for a perturbed
set-valued optimization problem. Inspired by Ha's work and using the
Gerstewitz's function (see, for example, [18-20]), the author [41]
obtained an improvement of Ha's version of EVP by relaxing several
assumptions. In the above Ha's and Qiu's versions, the perturbation
is given by a nonzero element $k_0$ of the ordering cone multiplied
by the distance function $d(\cdot, \cdot)$, i.e., its form is as
$d(\cdot,\cdot) k_0$ (disregarding a constant coefficient); and the
objective functions are set-valued maps. Bednarczuk and Zagrodny [7]
proved a vectorial EVP for a sequentially lower monotone
vector-valued map (which is called a monotonically semicontinuous
map in [7]), where the perturbation is given by a convex subset $H$
of the ordering cone multiplied by the distance function $d(\cdot,
\cdot)$, i.e., its form is as $d(\cdot, \cdot) H$. This generalizes
the case where directions of the perturbations are singletons $k_0$.
More generally, Guti\'{e}rrez, Jim\'{e}nez and Novo [23] introduced
a set-valued metric, which takes values in the set family of all
subsets of the ordering cone and satisfies the triangle inequality.
By using it they gave an original approach to extending the
scalar-valued EVP to a vector-valued map, where the perturbation
contains a set-valued metric. They also deduced several special
versions of EVP involving approximate solutions for vector
optimization problems and discussed their interesting applications
in optimization. In the above EVPs given by Bednarczuk and Zagrodny
[7] and by Guti\'{e}rrez, Jim\'{e}nez and Novo [23], the objective
maps are all a vector-valued (single-valued) map and the
perturbations  contain a convex subset of the ordering cone and a
set-valued metric with values in the ordering cone, respectively.

Very recently, Liu and Ng [33], Tammer and Z$\breve{a}$linescu [48]
and Flores-Ba\'{z}an, Guti\'{e}rrez and Novo [17] further considered
more general versions of EVP, where not only the objective map is a
set-valued map, but also the perturbation is a set-valued map, even
a family of set-valued maps satisfying certain property. In
particular, Liu and Ng [33] established several set-valued EVPs,
where the objective map is a set-valued map and the perturbation is
as the form $\gamma d(\cdot, \cdot) H$ or $\gamma^{\prime} d(\cdot,
\cdot) H,\ \gamma^{\prime}\in (0, \gamma)$, where $\gamma >0$ is a
constant, $d(\cdot, \cdot)$ is the metric  on the domain space and
$H$ is a closed convex subset of the ordering cone. Using the
obtained EVPs, they provided some sufficient conditions ensuring the
existence of error bounds for inequality systems. Tammer and
Z$\breve{a}$linescu [48] presented new minimal point theorems in
product spaces and the corresponding set-valued EVPs. As special
cases, they derived many of the previous EVPs and their extensions,
for example,  extensions of EVPs of Isac-Tammer's (see [28]) and
Ha's versions (see [24]). Through an extension of Br\'{e}zis-Browder
principle, Flores-Ba\'{z}an, Guti\'{e}rrez and Novo [17] established
a general strong minimal point existence theorem on quasi ordered
spaces and deduced several very general set-valued EVPs, where the
objective map is a set-valued map and the perturbation even involves
a family of set-valued maps satisfying ``triangle inequality"
property. As we have seen, these general set-valued EVPs extend and
improve the previous EVPs and imply many new interesting  results.

On the other hand, Bao and Mordukhovich (see [4,5]) proposed the
limiting monotonicity condition on objective maps and established
some enhanced versions of EVP for Pareto minimizers of set-valued
maps. By using minimal element theorems for product orders in
locally convex spaces, Khanh and Quy [31] generalized and improved
the above enhanced  versions of EVP. Particularly, they extended the
direction of the perturbation from a single positive vector to a
convex subset of the positive cone and removed the assumption in [4,
5] that the objective map is level closed.

 In this paper, we first establish a pre-order principle, which consists of a pre-order set $(X, \preceq)$
 and an extended real-valued function $\eta$ which is monotone with respect to $\preceq$.
 The pre-order principle states that there exists a strong minimal point dominated by any given point
 provided that the monotone function $\eta$  satisfies three general conditions.   From the pre-order principle,
  we obtain a very general set-valued EVP, where the objective function is a set-valued
map taking values in a quasi ordered linear space
 and the perturbation  contains a family of set-valued maps satisfying certain property.
Our assumption is accurate and weaker than ones appeared in the
previous EVPs. And our proof is clear and concise. The key to the
proof is to distinguish two different points by scalarizations. From
the general EVP, we can deduce all of the above mentioned set-valued
EVPs, their improvements and some new versions. In particular, our
pre-order principle also implies generalizations of Khanh and Quy's
minimal element theorems for product orders and hence we obtain
several versions of EVP for Pareto minimizers, which generalize and
improve the corresponding results of Bao and Mordukhovich  ([4, 5])
and of Khanh and Quy ([31]).

 The structure of this paper is as follows. In Section 2, we
 establish a pre-order principle. In
 Section 3, we give a general set-valued EVP and deduce a
 number of corollaries. In Section 4, we discuss set-valued EVPs,
 where perturbations contain a convex subset of the ordering cone.
 Moreover, we  give several set-valued EVPs for approximately
 efficient solutions.
In Section 5, we discuss minimal points for product orders and
present several versions of EVP for Pareto minimizers.\\

\section*{ \large\bf 2.  A pre-order  principle }

\hspace*{\parindent} Let $X$ be a nonempty set. As in [17], a binary
relation $\preceq$ on $X$ is called a pre-order if it satisfies the
transitive property; a quasi order if it satisfies the reflexive and
transitive properties; a partial order if it satisfies the
antisymmetric, reflexive and transitive properties. Let $(X,
\preceq)$ be a pre-order set. An extended real-valued function
$\eta:\, (X, \preceq) \rightarrow R\cup\{\pm \infty\}$ is called
monotone with respect to $\preceq$ if for any $x_1,\, x_2\in X$,
$$x_1\preceq x_2\ \ \Longrightarrow\ \ \eta(x_1)\leq \eta(x_2).$$
For any given $x_0\in X$, denote $S(x_0)$ the set $\{x\in X:\,
x\preceq x_0\}$. First we give a pre-order principle as follows.\\

{\bf Theorem 2.1.} \ {\sl Let $(X, \preceq)$ be a pre-order set,
$x_0\in X$ such that $S(x_0)\not=\emptyset$ and $\eta:\, (X,
\preceq) \rightarrow R\cup\{\pm\infty\}$ be an extended real-valued
function which is monotone with respect to $\preceq$.

Suppose that the following conditions are satisfied:

{\rm (A)}  \  $-\infty < \inf\{\eta(x):\, x\in S(x_0)\} <+\infty$.

{\rm (B)}  \  For any $x\in S(x_0)$ with $-\infty<\eta(x)<+\infty$
and $x^{\prime}\in S(x)\backslash\{x\}$, one has
$\eta(x)>\eta(x^{\prime})$.

{\rm (C)} \ For any sequence $(x_n)\subset S(x_0)$ with $x_n\in
S(x_{n-1}),\ \forall n$, such that  $\eta(x_n) -\inf\{\eta(x):\,
x\in S(x_{n-1})\} \rightarrow 0$ \  $(n\rightarrow \infty)$, there
exists $u\in X$ such that $u\in S(x_n),\ \forall n$.

Then there exists $\hat{x}\in X$ such that

{\rm (a)} \ $\hat{x}\in S(x_0)$;

{\rm (b)} \ $S(\hat{x})\subset \{\hat{x}\}$.}\\

{\bf Proof.}\  For brevity, we denote $\inf\{\eta(x):\, x\in
S(x_0)\}$ by $\inf \eta\circ S(X_0)$. By (A), we have
$$-\infty <\inf \eta\circ S(x_0)<+\infty. \eqno{(2.1)}$$
So, there exists $x_1\in S(x_0)$ such that
$$\eta(x_1) < \inf \eta\circ S(x_0) +\frac{1}{2}. \eqno{(2.2)}$$
By the transitive property of $\preceq$, we have
$$ S(x_1)\subset S(x_0). \eqno{(2.3)}$$
If $S(x_1)\subset \{x_1\}$, then we may take $\hat{x}:=x_1$ and
clearly $\hat{x}$ satisfies (a) and (b). If not, by (2.1), (2.2) and
(2.3) we conclude that
$$-\infty < \inf \eta\circ S(x_1) <+\infty.$$
So, there exists $x_2\in S(x_1)$ such that
$$\eta(x_2) < \inf\eta\circ S(x_1) +\frac{1}{2^2}.$$
In general, if $x_{n-1}\in X$ has been chosen, we may choose $x_n\in
S(x_{n-1})$ such that
$$\eta(x_n) < \inf\eta\circ S(x_{n-1}) +\frac{1}{2^n}.$$
If there exists $n$ such that $S(x_n)\subset \{x_n\}$, then we may
take $\hat{x} :=x_n$ and clearly $\hat{x}$ satisfies (a) and (b). If
not, we can obtain a sequence $(x_n)\subset S(x_0)$ with $x_n\in
S(x_{n-1}),\ \forall n$, such that
$$\eta(x_n) < \inf \eta\circ S(x_{n-1}) +\frac{1}{2^n},\ \forall n.
\eqno{(2.4)}$$ Obviously, $\eta(x_n) -\inf\eta\circ S(x_{n-1})
\rightarrow 0$ when $n\rightarrow \infty$. By (C), there exists
$\hat{x}\in X$ such that $$\hat{x}\in S(x_n), \ \ \forall n.
\eqno{(2.5)}$$ Obviously, $\hat{x}\in S(x_0)$, i.e., $\hat{x}$
satisfies (a). Next we show that $\hat{x}$ satisfies (b), i.e.,
$S(\hat{x})\subset\{\hat{x}\}$. If it is not, there exists
$\bar{x}\in S(\hat{x})$ and $\bar{x}\not=\hat{x}$. By (B),
$$\eta(\hat{x}) >\eta(\bar{x}). \eqno{(2.6)}$$
On the other hand, by $\bar{x}\in S(\hat{x})$ and (2.5) we have
$$\bar{x}\in S(x_n) \ \ \forall n.\eqno{(2.7)}$$
Since $\eta$ is monotone with respect to $\preceq$, by (2.5), (2.4)
and (2.7) we have
\begin{eqnarray*}
\eta(\hat{x})\leq\eta(x_n) \,&<&\,\inf\eta\circ S(x_{n-1})
+\frac{1}{2^n}\\
&\leq&\,\eta(\bar{x}) +\frac{1}{2^n},\ \forall n.
\end{eqnarray*}
Letting $n\rightarrow\infty$, we have
$\eta(\hat{x}) \leq\eta(\bar{x})$, which contradicts (2.6).  \hfill\framebox[2mm]{}\\

{\bf Remark 2.1.} \   The pre-order principle given by Theorem 2.1
consists of a pre-order set $(X, \preceq)$ and a monotone extended
real-valued function $\eta$ on $(X, \preceq)$. It states that there
exists a strong minimal point dominated by any given point provided
that the monotone function $\eta$  satisfies three general
conditions (A), (B) and (C). First, condition (A) is fundamental and
a starting point for constructing recurrently a decreasing sequence
$(x_n)$ with $\eta(x_n) -\inf \eta\circ S(x_{n-1}) \rightarrow 0\
(n\rightarrow\infty)$. Now that we have such a sequence $(x_n)$,
condition (C) plays a key role, which says that there exists
$\hat{x}\in X$ such that $\hat{x}\in S(x_n),\ \forall n$.
Particularly, $\hat{x}\in S(x_0)$ and  conclusion (a) holds. The
role of condition (B) is to distinguish points $x$ and non-$x$ in
$S(x)$.  Condition (B) together with the transitivity of $\preceq$
and the condition that $\eta(x_n) -\inf \eta\circ S(x_{n-1})
\rightarrow 0$  ensures that there is no $x^{\prime}\not=\hat{x}$
such that $x^{\prime}\in S(\hat{x})$. That is, conclusion (b) holds.
We realize that although the proofs of various versions of EVP may
be different, but their outlines are all similar to the above
process. We shall see that the pre-order principle indeed includes
many versions of EVP and their improvements. It should be noted also
that in [49, 50] various kinds of ordering principles were
established and many important applications were given. Our
pre-order principle is different from
them. It is specially made for deriving EVPs.\\

\section*{ \large\bf 3.  A general  set-valued EVP  and its corollaries }

\hspace*{\parindent} Let  $Y$ be a real linear space. If $A,\,
B\subset Y$ and $\alpha\in R$, the sets $A+B$ and $\alpha\,A$ are
defined as follows: $$A+B:=\{z\in Y:\, \exists x\in A,\, \exists
y\in B\ {\rm such\ that}\  z=x+y\},$$  $$\alpha A:=\{z\in Y:\,
z=\alpha x,\, x\in A\}.$$ A nonempty subset $D$ of $Y$ is called a
cone if $\alpha D\subset D$ for any $\alpha\geq 0$. And $D$ is
called a convex cone if $D+D\subset D$ and $\alpha D\subset D$ for
any $\alpha\geq 0$. A convex cone $D$ can specify a quasi order on
$Y$ as follows:
$$y_1,\, y_2\in Y,\ \ y_1\leq_D y_2\ \ \ \Longleftrightarrow\ \ \ \ y_1-y_2\in
-D.$$ In this case, $D$ is also called the ordering cone or positive
cone. We always assume that $D$ is nontrivial, i.e., $D\not=\{0\}$
and $D\not=Y$. An extended real function $\xi:\, Y\rightarrow
R\cup\{\pm\infty\}$ is said to be $D$-monotone if
$\xi(y_1)\leq\xi(y_2)$ whenever $y_1\leq_D y_2$. For any nonempty
subset $M$ of $Y$, we put $\inf \xi\circ M\,=\,\inf\{\xi(y):\, y\in
M\}$. If $\inf \xi\circ M
>-\infty$, we say that $\xi$ is lower bounded on $M$. For any given
$y\in Y$, sometimes we denote $\xi(y)$ by $\xi\circ y$. A family of
set-valued maps $F_{\lambda}:\, X\times X\rightarrow
2^D\backslash\{\emptyset\},\ \lambda\in\Lambda,$ is said to satisfy
the ``triangle inequality" property (briefly, denoted by property
TI, see [17]) if for each $x_i\in X,\ i=1,2,3,$ and
$\lambda\in\Lambda$ there exist $\mu,\nu\in\Lambda$ such that
$$F_{\mu}(x_1, x_2) + F_{\nu}(x_2, x_3)\,\subset\, F_{\lambda}(x_1,
x_3) +D.$$ Let $X$ be a nonempty set and let $f:\, X\rightarrow
2^Y\backslash\{\emptyset\}$ be a set-valued map. For any nonempty
set $A\subset X$, we put $f(A):=\cup\{f(x):\, x\in A\}$. For any
$x_1,\, x_2\in X$, define $x_2\preceq x_1$ iff
$$f(x_1)\subset f(x_2) + F_{\lambda}(x_2, x_1) + D,\ \
\forall\lambda\in\Lambda.$$\\

{\bf Lemma 3.1.}\ {\sl ``$\preceq$" is a pre-order on $X$, i.e., it
is a binary relation satisfying transitive property.}\\

{\bf Proof.}\  Let $x_2\preceq x_1$ and $x_3\preceq x_2$. We show
below that  $x_3\preceq x_1$. By the definition of $\preceq$, we
have
$$f(x_1)\subset f(x_2) + F_{\lambda} (x_2, x_1) +D,\ \
\forall\lambda\in\Lambda; \eqno{(3.1)}$$ and
$$ f(x_2)\subset f(x_3) + F_{\lambda}(x_3, x_2) +D, \ \
\forall\lambda\in\Lambda.\eqno{(3.2)}$$ For the above $x_1,\, x_2,\,
x_3\in X$ and any given $\lambda\in\Lambda$, there exists
$\mu,\,\nu\in\Lambda$ such that
$$F_{\mu}(x_3, x_2) + F_{\nu}(x_2, x_1)\subset F_{\lambda}(x_3, x_1)
+D.\eqno{(3.3)}$$ By (3.1),
$$f(x_1) \subset f(x_2) + F_{\nu}(x_2, x_1) +D. \eqno{(3.4)}$$
By (3.2),
$$f(x_2)\subset f(x_3) + F_{\mu}(x_3, x_2) +D.\eqno{(3.5)}$$
Combining (3.4), (3.5) and (3.3), we have
\begin{eqnarray*}
f(x_1)\,&\subset&\, f(x_2) +F_{\nu}(x_2, x_1) +D\\
&\subset&\, f(x_3) + F_{\mu}(x_3, x_2) +D + F_{\nu}(x_2, x_1) +D\\
&\subset&\, f(x_3) + F_{\lambda}(x_3, x_1) +D+D+D\\
&=&\, f(x_3) +F_{\lambda}(x_3, x_1) +D.
\end{eqnarray*}
Since $\lambda$ is arbitrary, we conclude that $x_3\preceq x_1$. \hfill\framebox[2mm]{}\\

{\bf Theorem 3.1.} \ {\sl Let $X$ be a nonempty set, $Y$ be a real
linear space, $D\subset Y$ be a convex cone specifying a quasi order
$\leq_D$ on $Y$, $f:\, X\rightarrow 2^Y\backslash\{\emptyset\}$ be a
set-valued map and $F_{\lambda}:\, X\times X\rightarrow
2^D\backslash\{\emptyset\},\ \lambda\in \Lambda$, be a family of
set-valued maps satisfying the property TI. Let $x_0\in X$ such that
$$S(x_0):=\{x\in X:\, f(x_0)\subset f(x) + F_{\lambda}(x, x_0) +D,\
\forall\lambda\in\Lambda\}\,\not=\,\emptyset.$$ Suppose that there
exists a $D$-monotone extended real function $\xi:\, Y\rightarrow
R\cup\{\pm\infty\}$ satisfying the following assumptions:

{\rm (D)} \  $-\infty  <\inf \xi\circ f(S(x_0))<+\infty$.

{\rm (E)} \   For any $x\in S(x_0)$ with $-\infty< \inf \xi\circ
f(x)<+\infty$ and for any $x^{\prime}\in S(x)\backslash\{x\}$, one
has $\inf\xi\circ f(x) >\inf\xi\circ f(x^{\prime})$.

{\rm (F)} \   For any  sequence $(x_n)\subset S(x_0)$ with $x_n\in
S(x_{n-1}),\ \forall n,$ such that  $\inf \xi\circ f(x_n)
-\inf\xi\circ f(S(x_{n-1})) \rightarrow 0,\ n\rightarrow\infty$,
there exists $u\in X$ such that $u\in S(x_n),\ \forall n.$

Then there exists $\hat{x}\in X$ such that

{\rm (a)} $f(x_0)\subset f(\hat{x}) +F_{\lambda}(\hat{x}, x_0) +D,\
\forall\lambda\in\Lambda$;

{\rm (b)} $\forall x\in X\backslash\{\hat{x}\}, \,
\exists\lambda\in\Lambda$ \ such \ that\ $f(\hat{x})\not\subset f(x)
+F_{\lambda}(x,\hat{x}) +D.$ }\\

{\bf Proof.}\  By Lemma 3.1, we can define a pre-order $\preceq$ on
$X$ as follows: for any $x_1,\, x_2\in X$,
$$x_2\preceq x_1 \ \ \  {\rm iff}\ \ \  f(x_1)\subset f(x_2) +
F_{\lambda}(x_2, x_1) +D,\ \forall\lambda\in\Lambda.$$ Thus,
$S(x_0)=\{x\in X:\, x\preceq x_0\}$. Define an extended real-valued
function  $\eta:\, (X, \preceq) \rightarrow R\cup\{\pm\infty\}$ as
follows
$$\eta(x):=\inf \xi\circ f(x),\ \ \forall x\in X.$$
Let $x^{\prime}\preceq x$. Then
$$f(x) \subset f(x^{\prime}) + F_{\lambda}(x^{\prime}, x) +D,\ \
\forall\lambda\in\Lambda.$$ For any $y\in f(x)$, there exists
$y^{\prime}\in f(x^{\prime}),\ q_{\lambda}(x^{\prime}, x) \in
F_{\lambda}(x^{\prime}, x),\ d_{\lambda, x^{\prime}, x}\in D$ such
that
$$y= y^{\prime} + q_{\lambda}(x^{\prime}, x) + d_{\lambda,
x^{\prime}, x} .$$ Since $$y-y^{\prime} = q_{\lambda}(x^{\prime}, x)
+d_{\lambda, x^{\prime}, x} \in D,$$ we have
$$\xi(y) \geq \xi(y^{\prime}) \geq \inf\xi\circ f(x^{\prime}).$$
As $y\in f(x)$ is arbitrary, we have
$$\inf \xi\circ f(x) \geq \inf \xi\circ f(x^{\prime}),\ \  {\rm
i.e.,}\  \  \eta(x)\geq\eta(x^{\prime}).$$ Thus, $\eta$ is monotone
with respect to $\preceq$. It is easy to see that assumptions (D),
(E) and (F) are exactly assumptions (A), (B) and (C) in Theorem 2.1.
Now, applying Theorem 2.1, we know that there exists $\hat{x}\in X$
such that $\hat{x}\in S(x_0)$ and $S(\hat{x})\subset \{\hat{x}\}$.
This means that
$$f(x_0)\subset f(\hat{x}) +F_{\lambda}(\hat{x}, x_0)
+D,\ \forall\lambda\in\Lambda$$ and
 $$\forall x\in X\backslash\{\hat{x}\}, \, \exists\lambda\in\Lambda
\ {\rm such \ that}\ f(\hat{x})\not\subset f(x)
+F_{\lambda}(x,\hat{x}) +D.$$ That is, (a) and (b) are satisfied. \hfill\framebox[2mm]{}\\

For a real linear space $Y$, denote the algebraic dual of $Y$ by
$Y^{\#}$ and  denote the positive polar cone of  $D$ in $Y^{\#}$ by
$D^{+\#}$, i.e., $D^{+\#}=\{l\in Y^{\#}:\, l(d)\geq 0,\,\forall d\in
D\}$. Obviously, every $\xi\in D^{+\#}$ is a $D$-monotone real
function. Hence, the $\xi$ in Theorem 3.1 may be an element of
$D^{+\#}\backslash \{0\}$. In this case, assumptions (D) and (E)
become more concise. And the expression of assumption (F) is the same as in Theorem 3.1.\\

{\bf Theorem  3.1$^{\prime}$.} \ {\sl  let $X,\,Y,\, D,\, f,\,
F_{\lambda},\,\lambda\in\Lambda$, and $x_0\in X$ be the same as in
Theorem 3.1. Suppose that there exists $\xi\in
D^{+\#}\backslash\{0\}$ satisfying the following assumptions:

{\rm (D)} \  $\xi$ is lower bounded on $f(S(x_0))$, i.e., $-\infty <
\inf \xi\circ f(S(x_0))$.

{\rm (E)} \  For any $x\in S(x_0)$ and any $x^{\prime}\in
S(x)\backslash\{x\}$, one has $\inf \xi\circ f(x) > \inf \xi\circ
f(x^{\prime})$.

{\rm (F)} \ See Theorem 3.1.

Then the result of Theorem 3.1 remains true.}\\

From Theorem 3.1$^{\prime}$, we can deduce [17, Theorem 6.5], a
general version of set-valued EVP, which extends EVPs in [7, 23].\\

{\bf Corollary 3.1} (see [17, Theorem 6.5]).\ {\sl The result of
Theorem 3.1$^{\prime}$ remains true if assumption {\rm (E)} is
replaced by the following assumption

{\rm (E$_1$)} \   For any $x,\, x^{\prime}\in S(x_0)$ with $x\not=
x^{\prime}$, there exists $\lambda_0\in\Lambda$ such that
$\inf\{\xi\circ q_{\lambda_0}:\,
q_{\lambda_0}\in F_{\lambda_0}(x^{\prime}, x)\} > 0$.}\\

{\bf Proof.}\ By Theorem  3.1$^{\prime}$ we only need to prove that
 (E$_1$) $\Rightarrow$ (E). Let $x\in S(x_0)$ and
$x^{\prime}\in S(x)\backslash\{x\}$. We shall show that $\inf
\xi\circ f(x)\,>\,\inf\xi\circ f(x^{\prime})$. Obviously,
$$x^{\prime}\not=x\ \ {\rm and}\ \ f(x)\subset f(x^{\prime}) +
F_{\lambda}(x^{\prime}, x) +D,\ \ \forall\lambda\in\Lambda.$$ By
(E$_1$), there exists $\lambda_0\in\Lambda$ such that
$$\eta:=\inf\{\xi\circ q_{\lambda_0}:\, q_{\lambda_0}\in
F_{\lambda_0}(x^{\prime}, x)\}\,>\, 0.$$ Clearly, there exists $y\in
f(x)$ such that
$$\xi\circ y\,<\,\inf\xi\circ f(x) +\frac{1}{2}\eta.\eqno{(3.6)}$$
As $y\in f(x)\subset f(x^{\prime}) +F_{\lambda_0}(x^{\prime}, x)
+D$, there exists $y^{\prime}\in f(x^{\prime})$, $q_{\lambda_0}\in
F_{\lambda_0}(x^{\prime}, x)$ and $d\in D$ such that $$y\,=\,
y^{\prime} +q_{\lambda_0} +d.$$ Since $\xi\in
D^{+\#}\backslash\{0\}$, we have
$$\xi\circ y \,=\,\xi\circ y^{\prime} +\xi\circ q_{\lambda_0}
+\xi\circ d\,\geq\,\xi\circ y^{\prime} +\eta.\eqno{(3.7)}$$
Combining (3.6) and (3.7), we have
$$\inf\xi\circ f(x) +\frac{1}{2}\eta \,>\,\xi\circ y\,\geq\,\xi\circ
y^{\prime} +\eta\,\geq\,\inf\xi\circ f(x^{\prime}) +\eta.$$ Thus,
$$\inf \xi\circ f(x)\,>\,\inf\xi\circ
f(x^{\prime})+\frac{1}{2}\eta\,>\,\inf\xi\circ f(x^{\prime}).$$\hfill\framebox[2mm]{}\\

If $\xi$ attains its infimum on $f(x)$ for any $x\in S(x_0)$, then
assumption (E$_1$) in Corollary 3.1 can be weakened.\\

{\bf Corollary 3.2.} \ {\sl Suppose that $\xi$ in Theorem
3.1$^{\prime}$ attains its infimum on $f(x)$ for any $x\in S(x_0)$.
Then the result of Theorem 3.1$^{\prime}$ remains true if assumption
{\rm (E)} is replaced by the following assumption

{\rm (E$_2$)} \  For any $x,\,x^{\prime}\in S(x_0)$ with
$x\not=x^{\prime}$, there exists $\lambda_0\in\Lambda$ such that
$\xi\circ q_{\lambda_0}>0,\ \forall q_{\lambda_0}\in
F_{\lambda_0}(x^{\prime}, x)$.}\\

 {\bf Proof.}\  By Theorem 3.1$^{\prime}$ we only need to prove that (E$_2$)
$\Rightarrow$ (E). Let $x\in S(x_0)$ and $x^{\prime}\in
S(x)\backslash\{x\}$. Then $x^{\prime}\not=x$  and  $$f(x)\subset
f(x^{\prime}) + F_{\lambda}(x^{\prime}, x) +D,\ \
\forall\lambda\in\Lambda.\eqno{(3.8)}$$ By assumption (E$_2$), there
exists $$y_x\in f(x)\eqno{(3.9)}$$ such that
$$\xi\circ y_x\,=\,\inf\xi\circ f(x). \eqno{(3.10)}$$
Since $x\not= x^{\prime}$, by (E$_2$) there exists
$\lambda_0\in\Lambda$ such that
$$\xi\circ q_{\lambda_0} >0,\ \ \  \forall q_{\lambda_0}\in
F_{\lambda_0}(x^{\prime}, x).$$ By (3.9) and (3.8),
$$y_x\,\in\, f(x^{\prime}) +F_{\lambda_0}(x^{\prime}, x) +D.$$
Hence, there exists $y^{\prime}\in f(x^{\prime}),\ q_{\lambda_0}\in
F_{\lambda_0}(x^{\prime}, x)$ and $d\in D$ such that
$$y_x\,=\, y^{\prime} +q_{\lambda_0} +d.$$
As $\xi\in D^{+\#}$ and $\xi\circ q_{\lambda_0} >0$, we have
\begin{eqnarray*}
\xi\circ y_x \,&=&\,\xi\circ y^{\prime} +\xi\circ q_{\lambda_0}
+\xi\circ d \\
&\geq&\, \xi\circ y^{\prime} +\xi\circ q_{\lambda_0}\\
 &>&\,\xi\circ
y^{\prime}\\
&\geq&\,\inf\xi\circ f(x^{\prime}).
\end{eqnarray*}
Combining this with (3.10), we conclude that $\inf \xi\circ f(x)
>\inf
\xi\circ f(x^{\prime})$ and hence (E) holds. \hfill\framebox[2mm]{}\\

In particular, if $f(x)$ is a singleton for any $x\in X$, i.e.,
$f:\, X\rightarrow Y$ is a vector-valued map, then from Corollary
3.2 we can obtain [42, Theorem 3.15]. If $Y$ is a separated
topological vector space, we denote $Y^*$ the topological dual of
$Y$ and denote $D^+$ the positive polar cone of $D$ in $Y^*$, i.e.,
$D^+:=\{l\in Y^*:\, l(d)\geq 0,\,\forall d\in D\}$. It is possible
that $Y^*=\{0\}$, i.e., there is no non-trivial continuous linear
functional on $Y$. However, if $Y$ is a locally convex separated
topological vector space (briefly, denoted by a locally convex
space), then $Y^*$ is large enough so that $Y^*$ can separate points
in $Y$. Suppose that $K\subset Y$ is a weakly countably compact set
(particularly, a weakly compact set) in the locally convex space
$Y$. Then for any $l\in Y^*$, $l$ attains its infimum
on $K$. Thus, from Corollary 3.2,  we have the following.\\

{\bf Corollary 3.3.} \ {\sl In Theorem 3.1$^{\prime}$, we further
assume that $Y$ is a locally convex space and the set-valued map
$f:\, X\rightarrow 2^Y\backslash\{\emptyset\}$ satisfies that $f(x)$
is a weakly countably compact set (particularly, a weakly compact
set) in $Y$ for any $x\in S(x_0)$. Then the result of Theorem
3.1$^{\prime}$ remains true if assumption {\rm (E)} is replaced by
{\rm (E$_2$)}}\\

{\bf Definition 3.1} (see [42, Definition 3.4]).  \   Let $X$ be a
topological space and let $S(\cdot): X\rightarrow
2^X\backslash\{\emptyset\}$ be a set-valued map. The set-valued map
$S(\cdot)$ is said to be dynamically closed at $x\in X$ if
$(x_n)\subset S(x),\ S(x_{n+1})\subset S(x_n)\subset S(x)$ for all
$n$ and $x_n\rightarrow \bar{x}$ then $\bar{x}\in S(x)$. In this
case, we also say that $S(x)$ is dynamically closed. Moreover, let
$(X, d)$ be a metric space and $x\in X$. Then $(X, d)$ is said to be
$S(x)$-dynamically complete if every Cauchy
sequence $(x_n)\subset S(x)$ such that $S(x_{n+1})\subset S(x_n)\subset S(x)$ for all $n$, is convergent in $X$.\\

The following corollary generalizes [48,  Theorems 4.1 and
6.1].\\

{\bf Corollary 3.4.} \ {\sl Let $(X, d)$ be a complete metric space,
$Y$ be a real linear space, $D\subset Y$ be a convex cone specifying
a quasi order $\leq_D$ on $Y$, $f:\, X\rightarrow
2^Y\backslash\{\emptyset\}$ be a set-valued map and $F_{\lambda}:\,
X\times X \rightarrow 2^D\backslash\{\emptyset\}$,
$\lambda\in\Lambda,$ be a family of set-valued maps satisfying the
property TI.  Let $x_0\in X$ such that $S(x_0):=\{x\in X:\,
f(x_0)\subset f(x) +F_{\lambda}(x, x_0) +D,\
\forall\lambda\in\Lambda\}\not=\emptyset$. Suppose that for any
$x\in S(x_0)$, $S(x)$ is dynamically closed and that there exists
$\xi\in D^{+\#}\backslash\{0\}$ satisfying the following
assumptions:

{\rm (D)} \ $\xi$ is lower bounded on $f(S(x_0))$.

{\rm (E$_3$)} \   There exists $\lambda_0\in\Lambda$ such that for
any $\delta >0$, $\inf\xi\circ F_{\lambda_0\delta}\,>\,0$, where
$F_{\lambda_0\delta}$ denotes the set $\cup_{d(x,
x^{\prime})\geq\delta}F_{\lambda_0}(x, x^{\prime})$.

Then there exists $\hat{x}\in X$ such that

{\rm (a)} $f(x_0)\,\subset\, f(\hat{x}) +F_{\lambda}(\hat{x}, x_0)
+D,\ \ \forall \lambda\in\Lambda;$

{\rm (b)} $\forall x\in X\backslash\{\hat{x}\},\
\exists\lambda\in\Lambda$ such that $f(\hat{x})\,\not\subset\, f(x)
+F_{\lambda}(x, \hat{x}) +D.$ }\\

 {\bf Proof.}\  Obviously, (E$_3$) $\Rightarrow$ (E$_1$). By
 Corollary 3.1, we only need to prove that assumption (F) in Theorem 3.1$^{\prime}$ is
 satisfied.  Let $(x_n)\subset S(x_0)$ such that $x_n\in S(x_{n-1})$
 and
 $$\inf \xi\circ f(x_n)\,-\, \inf\xi\circ f(S(x_{n-1}))
 \,\rightarrow\, 0,\ \ \ n\rightarrow \infty.$$
We may take a positive real sequence $(\epsilon_n)$ such that
$\epsilon_n \rightarrow 0 \ (n\rightarrow\infty)$ and such that
$$\inf \xi\circ f(x_n) \,<\, \inf \xi\circ f(S(x_{n-1}))
+\epsilon_n, \ \ \forall n.$$
For each $n$, choose $y_n\in f(x_n)$
such that
$$\xi\circ y_n\,<\,\inf\xi\circ f(S(x_{n-1}))
+\epsilon_n.\eqno{(3.11)}$$ We assert that $(x_n)$ is a Cauchy
sequence in $(X, d)$. If not, there exists $\delta >0$ such that for
any $k\in N$, there exists $n_k>k$ such that $$d(x_{n_k},
x_k)\,\geq\, \delta.\eqno{(3.12)}$$ For any $k\in N$, $x_{n_k}\in
S(x_{n_k-1})\subset S(x_k)\subset S(x_{k-1})$. Thus,
$$y_k\in f(x_k)\subset f(x_{n_k}) + F_{\lambda_0}(x_{n_k}, x_k)
+D.$$ Hence, there exists $y_{n_k}^{\prime}\in f(x_{n_k})$,
$q_{\lambda_0, n_k, k}\in F_{\lambda_0}(x_{n_k}, x_k)$ and $d_k\in
D$ such that
$$y_k\,=\, y_{n_k}^{\prime} +q_{\lambda_0, n_k, k}+d_k.$$
As $\xi\in D^{+\#}\backslash\{0\}$, we have
$$\xi\circ y_k\,=\,\xi\circ y_{n_k}^{\prime} +\xi\circ q_{\lambda_0,
n_k, k}+\xi\circ d_k\,\geq\,\xi\circ y_{n_k}^{\prime} +\xi\circ
q_{\lambda_0, n_k, k}.\eqno{(3.13)}$$ Remarking that
$y_{n_k}^{\prime}\in f(x_{n_k})\subset f(S(x_{k-1}))$,  we have
$$\xi\circ y_{n_k}^{\prime}\,\geq\,\inf\xi\circ
f(S(x_{k-1})).\eqno{(3.14)}$$ By (3.13), (3.14) and (3.11), we have
\begin{eqnarray*}
\xi\circ q_{\lambda_0, n_k, k}\,&\leq&\, \xi\circ y_k -\xi\circ
y_{n_k}^{\prime}\\
&\leq&\,\xi\circ y_k-\inf\xi\circ f(S(x_{k-1}))\\ &<&\,\epsilon_k.
\end{eqnarray*}
Letting $k\rightarrow\infty$, we have
$$\xi\circ q_{\lambda_0, n_k, k}\,\rightarrow\,0\ \
(k\rightarrow\infty).\eqno{(3.15)}$$ On the other hand, by (3.12),
every $q_{\lambda_0, n_k, k}\in F_{\lambda_0}(x_{n_k}, x_k)\subset
F_{\lambda_0\delta}.$ By (E$_3$),
$$\xi\circ q_{\lambda_0, n_k, k}\,\geq\,\inf\xi\circ
F_{\lambda_0\delta}\,>\,0,\ \ \forall k,$$ which contradicts (3.15).
Thus, we have shown that  $(x_n)$ is a Cauchy sequence. Since $(X,
d)$ is complete, there exists $u\in X$ such that $x_n\rightarrow u$.
By the assumption, $S(x_n)$ is dynamically closed  for any $n$.
Since $(x_{n+p})_{p\in N}\subset S(x_n),\  x_{n+p+1}\in S(x_{n+p})$
and $x_{n+p}\rightarrow u\ (p\rightarrow\infty)$, we have $u\in
S(x_n),\ \forall n.$ This means that assumption (F) is satisfied. \hfill\framebox[2mm]{}\\

Let $X$ be a metric space, $Y$ be a locally convex space and
$D\subset Y$ be a closed convex cone. As in [24], a set-valued map
$f:\, X\rightarrow 2^Y\backslash\{\emptyset\}$ is said to be
$D$-lower semicontinuous (briefly, $D$-l.s.c.) on $X$ if for any
$y\in Y$, the set $\{x\in X:\, f(x)\cap(y-D)\not=\emptyset\}$ is
closed. And $f:\, X\rightarrow 2^Y\backslash\{\emptyset\}$ is said
to have $D$-closed values if for any $x\in X$, $f(x)+D$ is closed.
$f(X)$ is said to be $D$-bounded if there exists a bounded set
$M\subset Y$ such that $f(X)\subset M+D$. In [48], $f(X)$ being
$D$-bounded is also called $f(X)$  being quasibounded from below. In
[41], a set-valued map $f:\, X\rightarrow
2^Y\backslash\{\emptyset\}$ is said to be $D$-sequentially lower
monotone (briefly, $D$-s.l.m.) if $f(x_n)\subset f(x_{n+1})+D,\
\forall n,$ and $x_n\rightarrow\bar{x}$ imply $f(x_n)\subset
f(\bar{x}) +D,\ \forall n$.  It is easy to show that a $D$-l.s.c.
set-valued map is $D$-s.l.m. But the converse is not true (see
[41]). From Corollary 3.4 we can deduce an improvement of Ha's
version of set-valued EVP (see [24, Theorem 3.1]) as follows.\\

{\bf Corollary 3.5.} \ {\sl Let $(X, d)$ be a complete metric space,
$Y$ be a locally convex space pre-ordered by a closed convex cone
$D$ and $k_0\in D\backslash -D$. Let $f:\, X\rightarrow
2^Y\backslash\{\emptyset\}$ be $D$-s.l.m., have $D$-closed values
and $f(X)$ be $D$-bounded. Suppose that $x_0\in X$ and $\epsilon >0$
such that $f(x_0)\not\subset f(x) +\epsilon k_0 +D,\ \forall x\in
X$. Then for any $\lambda>0$, there exists $\hat{x}\in X$ such that

{\rm (a)} $f(x_0)\subset f(\hat{x}) +(\epsilon/\lambda) d(\hat{x},
x_0) k_0 +D$;

{\rm (b)} $\forall x\in X\backslash\{\hat{x}\},\
f(\hat{x})\not\subset f(x) +(\epsilon/\lambda) d(x, \hat{x}) k_0
+D.$, i.e., $\hat{x}$ is a strict minimizer of the map $x\mapsto
f(x) + (\epsilon/\lambda) d(x, \hat{x})k_0$ {\rm (concerning a
strict minimizer of a map, see  [24])};

{\rm (c)} $d(x_0, \hat{x})\leq\lambda$.}\\

 {\bf Proof.}\  For any $x,\, x^{\prime}\in X$, define $F(x,
 x^{\prime}):=(\epsilon/\lambda)d(x, x^{\prime}) k_0$. Obviously,
 the family $\{F\}$ satisfies the property TI. Put $$S(x_0):=\{x\in
 X:\, f(x_0)\subset f(x)+\frac{\epsilon}{\lambda} d(x, x_0)k_0
 +D\}.$$ Then $x_0\in S(x_0)$ and $S(x_0)\not=\emptyset$.
Since $k_0\not\in -D$ and $D$ is closed, by the separation theorem,
there exists $\xi\in D^+\backslash\{0\}$ such that $\xi(k_0)=1$.
When $d(x, x^{\prime})\geq \delta$, we have
$$\xi(\frac{\epsilon}{\lambda} d(x, x^{\prime}) k_0)\,\geq\,
\frac{\epsilon}{\lambda}\cdot \delta \,>\, 0.$$ Hence
$$\inf\xi\circ F_{\delta}\,>\, 0,\ \ {\rm where}\ \ F_{\delta}=\cup_{d(x,
x^{\prime})\geq\delta}F(x, x^{\prime})=\{\frac{\epsilon}{\lambda}
d(x, x^{\prime}) k_0:\, d(x, x^{\prime})\geq\delta\}.$$ Remarking
that $f(X)$ is $D$-bounded and $\xi\in D^+\backslash\{0\}$, we see
that $\xi$ is lower bounded on $f(X)$ and lower bounded on
$f(S(x_0))$. Thus, we have verified that assumptions (D) and (E$_3$)
in Corollary 3.4 are satisfied.

It remains to show that for any $x\in S(x_0)$, $S(x)$ is dynamically
closed. Let $(x_n)\subset S(x),\ x_{n+1}\in S(x_n)\subset S(x),\
\forall n,$ and $x_n\rightarrow u$. We shall show that $u\in S(x)$.
Since $x_{n+1}\in S(x_n)$,
$$f(x_n)\,\subset\, f(x_{n+1}) +\frac{\epsilon}{\lambda} d(x_{n+1},
x_n) k_0 +D\,\subset\,f(x_{n+1})+D,\ \ \forall n.$$ Combining this
with $x_n\rightarrow u$ and using the assumption that $f$ is
$D$-s.l.m., we have
$$f(x_n)\,\subset\, f(u) +D,\ \ \forall n.\eqno{(3.16)}$$
By $x_n\in S(x)$ and (3.16), we have
\begin{eqnarray*}
f(x)\,&\subset&\, f(x_n) +\frac{\epsilon}{\lambda} d(x_n, x) k_0
+D\\
&\subset&\, f(u) +D+\frac{\epsilon}{\lambda} d(x_n, x) k_0 +D\\
&=&\, f(u)+\frac{\epsilon}{\lambda} d(x_n, x) k_0 +D.
\end{eqnarray*}
Thus,
$$f(x)-\frac{\epsilon}{\lambda} d(x_n, x)k_0\,\subset\, f(u)
+D.\eqno{(3.17)}$$ Since $(\epsilon/\lambda)d(x_n, x) k_0\rightarrow
(\epsilon/\lambda) d(u, x) k_0\ (n\rightarrow\infty)$ and $f(u) +D$
is closed, by (3.17) we have
$$f(x) -\frac{\epsilon}{\lambda} d(u, x) k_0\,\subset\, f(u) +D,$$
which means that $u\in S(x)$.

Now applying Corollary 3.4, there exists $\hat{x}\in X$ such that
(a) and (b) hold. It remains to show that (c) holds. If not, suppose
that $d(\hat{x}, x_0)>\lambda$. Then by (a), we have
\begin{eqnarray*}
f(x_0)\,&\subset&\, f(\hat{x}) +\frac{\epsilon}{\lambda} d(\hat{x},
x_0) k_0 +D\\
&=&\, f(\hat{x}) +(\frac{\epsilon}{\lambda} d(\hat{x},
x_0)-\epsilon) k_0 +\epsilon k_0 +D\\
&\subset&\, f(\hat{x}) +\epsilon k_0 +D.
\end{eqnarray*}
This contradicts the assumption that $f(x_0)\not\subset f(x)
+\epsilon k_0 +D,\ \forall x\in X$. \hfill\framebox[2mm]{}\\

In fact, the assumption that $f(X)$ is $D$-bounded in Corollary 3.5
can be replaced by a weaker assumption: there exists $\epsilon
>0$ such that $f(x_0)\not\subset f(X) +\epsilon k_0 +D$ (see [41]).
For this, we need the following nonlinear scalarization function.
The original version is due to Gerstewitz [18].

We present the concept in a general setting. Let $Y$ be a real
linear space and $A\subset Y$ be a nonempty set. We put (refer to
[1])
$${\rm vcl} (A):=\{y\in Y:\, \exists v\in Y, \exists \lambda_n\geq
0,\,\lambda_n\rightarrow 0\  {\rm such\ that}\ y+\lambda_n v\in A,\
\forall n\in N\}.$$ For any given $v_0\in Y$, put
$${\rm vcl}_{v_0}(A)\,=\,\{y\in Y:\,\exists \lambda_n\geq 0,\,
\lambda_n\rightarrow 0 \ {\rm such\ that}\ y+\lambda_n v_0\in A,\
\forall n\in N\}.$$ Obviously, $$A\subset {\rm vcl}_{v_0}(A)\subset
\bigcup_{v\in Y} {\rm vcl}_v(A) ={\rm vcl}(A).$$ Moreover, if $Y$ is
a topological vector space, then ${\rm vcl}(A)\subset {\rm cl}(A)$
and the inclusion is proper. If $A={\rm vcl}_{v_0}(A)$, then $A$ is
said  to be $v_0$-closed. If $A={\rm vcl}(A)$, then $A$ is said to
be vectorial closed. Obviously, $v_0$-closedness and vectorial
closedness are both strictly weaker than topological closedness.

Let $D\subset Y$ be a  convex cone specifying a quasi order $\leq_D$
on $Y$ and $k_0\in D\backslash -{\rm vcl}(D)$. For any $y\in Y$, if
there exists $t\in R$ such that $y\in tk_0-D$, then for any
$t^{\prime}>t$, $y\in t^{\prime} k_0-D$. We define a function
$\xi_{k_0}:\, Y\rightarrow R\cup\{+\infty\}$ as follows: if there
exists $t\in R$ such that $y\in tk_0-D$, then define $\xi_{k_0}(y)
=\inf\{t\in R:\, y\in tk_0-D\}$; or else define
$\xi_{k_0}(y)=+\infty$.  As $k_0\not\in -{\rm vcl} (D)$, we can show
that $\xi_{k_0}(y)\not=-\infty$. This function is called a
Gerstewitz's function. Concerning the details of such a function and
its properties, please refer to [18-20]. For brevity, we denote
$D+(0, +\infty)k_0$ by ${\rm vint}_{k_0}(D)$. We list several
properties of
$\xi_{k_0}$ as follows.\\

{\bf Lemma 3.2} (refer to [10, 21,  40,  41, 44]). \ {\sl Let $y\in
Y$ and $r\in R$. Then we have:

{\rm (i)} $\xi_{k_0}(y)<r\ \Leftrightarrow\ y\in r\,k_0-{\rm
vint}_{k_o}(D)$.

{\rm (ii)} $\xi_{k_0}(y)\leq r\ \Leftrightarrow\ y\in rk_0-{\rm
vcl}_{k_0}(D).$

{\rm (iii)} $\xi_{k_0}(y) = r\  \Leftrightarrow\ y\in rk_0-\left(
{\rm vcl}_{k_0}(D)\backslash {\rm vint}_{k_0}(D)\right)$.
 In particular, $\xi_{k_0}(k_0)=1$
and $\xi_{k_0}(0)=0$.

{\rm (iv)} $\xi_{k_0}(y)\geq r\ \Leftrightarrow\ y\not\in rk_0-{\rm
vint}_{k_0}(D)$.

 {\rm (v)} $\xi_{k_0}(y)>r\ \Leftrightarrow\
y\not\in rk_0 -{\rm vcl}_{k_0}(D)$.

Moreover, we have:

{\rm (vi)} $\xi_{k_0}(y_1+y_2)\leq \xi_{k_0}(y_1) +\xi_{k_0}(y_2),\
\forall y_1,\, y_2\in Y.$

{\rm (vii)} $\xi_{k_0}(y+\lambda k_0) = \xi_{k_0}(y) +\lambda,\
\forall y\in Y,\, \forall\lambda\in R$.

{\rm (viii)} $y_1\leq_D y_2\ \Rightarrow\ \xi_{k_0}(y_1)\leq
\xi_{k_0}(y_2).$}\\

In Corollary 3.5, we need to assume that $f:\, X\rightarrow
2^Y\backslash\{\emptyset\}$  has $D$-closed values, i.e., $f(x)+D$
is closed for all $x\in X$. Here, we introduce a weaker property: a
set-valued map $f$ is said to have $D$-$k_0$-closed values if $f(x)
+D$ is $k_0$-closed for all $x\in X$.\\

{\bf Corollary 3.6} (see [41, Theorem 3.1]).  \ {\sl Let $(X, d)$ be
a metric space, $Y$ be a locally convex space quasi ordered by a
convex cone $D$ and $k_0\in D\backslash -{\rm vcl}(D)$. Let $f:\,
X\rightarrow 2^Y\backslash\{\emptyset\}$ be $D$-s.l.m. and have
$D$-$k_0$-closed values. Suppose that $x_0\in X$ and $\epsilon >0$
such that $f(x_0)\not\subset f(X) +\epsilon k_0 +D$ and suppose that
$(X, d)$ is $S(x_0)$-dynamically complete, where $S(x_0)=\{x\in X:
f(x_0)\subset f(x) +(\epsilon/\lambda) d(x, x_0) k_0 +D\}$. Then for
any $\lambda>0$, there exists $\hat{x}\in X$ such that

{\rm (a)} $f(x_0)\subset f(\hat{x}) +(\epsilon/\lambda) d(\hat{x},
x_0) k_0 +D$;

{\rm (b)} $\forall x\in X\backslash\{\hat{x}\},\
f(\hat{x})\not\subset f(x) +(\epsilon/\lambda) d(x, \hat{x}) k_0
+D;$

{\rm (c)} $d(x_0, \hat{x})\leq\lambda$.
}\\

{\bf Proof.}\  We shall prove the result by using Theorem 3.1.
Define $F:\, X\times X\rightarrow 2^D\backslash\{\emptyset\}$ as
follows: $F(x, x^{\prime}):=(\epsilon/\lambda)d(x, x^{\prime}) k_0$.
Then the family $\{F\}$ satisfies the property TI. Put
$$S(x_0):=\{x\in X:\, f(x_0)\subset f(x)+(\epsilon/\lambda)d(x, x_0) k_0
+D\}.$$ Obviously, $x_0\in S(x_0)$ and $S(x_0)\not=\emptyset$.
 Since $f(x_0)\not\subset f(S(x_0)) +\epsilon k_0 +D$, there exists $y_0\in
f(x_0)$ such that $y_0\not\in f(S(x_0)) +\epsilon k_0 +D$, that is,
$$ (f(S(x_0))-y_0)\cap (-\epsilon k_0 -D)\,=\,\emptyset.\eqno{(3.18)}$$
By (3.18) and Lemma 3.2(iv), we have
$$\xi_{k_0}(y-y_0)\geq -\epsilon,\ \ \forall y\in f(S(x_0)).$$
Also, $\xi_{k_0}(y_0-y_0)=0$, where $y_0\in f(x_0)\subset
f(S(x_0))$. Put $$\xi(y)=\xi_{k_0}(y-y_0),\ \ \  \forall y\in Y.$$
Then $\xi$ is a $D$-monotone extended real function such that
$$-\infty<-\epsilon \leq \inf\xi\circ f(S(x_0))\leq 0<+\infty.$$
That is, assumption (D) in Theorem 3.1 is satisfied. For any $x\in
S(x_0)$, we remark that
$$\inf \xi\circ f(x)\leq\inf\xi\circ f(x_0)\leq
\xi(y_0)=0<+\infty.$$ For any $x^{\prime}\in S(x)\backslash\{x\}$,
we have
$$f(x)\,\subset\, f(x^{\prime})+\frac{\epsilon}{\lambda}
d(x^{\prime}, x) k_0 +D. \eqno{(3.19)}$$  Choose $y\in f(x)$ such
that
$$\xi\circ y\,<\, \inf\xi\circ f(x) +\frac{\epsilon}{2\lambda}
d(x^{\prime}, x).\eqno{(3.20)}$$ By (3.19),
$$y-y_0\,\in\, f(x^{\prime})-y_0 +\frac{\epsilon}{\lambda}
d(x^{\prime}, x) k_0 +D.$$
Thus, there exists $y^{\prime}\in
f(x^{\prime})$ and $d\in D$ such that
$$y-y_0\,=\, y^{\prime}-y_0 +\frac{\epsilon}{\lambda} d(x^{\prime},
x) k_0 +d.$$ By Lemma 3.2, we have
$$\xi_{k_0}(y-y_0)\,\geq\, \xi_{k_0}\left(y^{\prime}-y_0+\frac{\epsilon}{\lambda}d(x^{\prime}, x)k_0\right)\,=\,
\xi_{k_0} (y^{\prime} -y_0)
+\frac{\epsilon}{\lambda} d(x^{\prime}, x),$$ that is,
$$\xi\circ y\,\geq\, \xi\circ y^{\prime} +\frac{\epsilon}{\lambda}
d(x^{\prime}, x).\eqno{(3.21)}$$ Combining (3.20) and (3.21), we
have
$$\inf \xi\circ f(x) +\frac{\epsilon}{2\lambda} d(x^{\prime} , x)
>\xi\circ y \geq \xi\circ y^{\prime} +\frac{\epsilon}{\lambda}
d(x^{\prime}, x)\geq \inf\circ f(x^{\prime})
+\frac{\epsilon}{\lambda} d(x^{\prime}, x).$$ From this,
$$\inf\xi\circ f(x) >\inf\xi\circ f(x^{\prime})
+\frac{\epsilon}{2\lambda}  d(x^{\prime}, x) >\inf\xi\circ
f(x^{\prime}).$$ That is, assumption (E) in Theorem 3.1 is
satisfied.

Finally, we show that assumption (F) in Theorem 3.1 is satisfied.
Let a sequence $(x_n)\subset S(x_0)$ with $x_n\in S(x_{n-1}),\
\forall n,$ such that $$\inf\,\xi\circ f(x_n) -\inf\,\xi\circ
f(S(x_{n-1}))\,\rightarrow\, 0,\ n\rightarrow\infty.$$ Put
$$\epsilon_n:=\,\inf\,\xi\circ f(x_n) - \inf\,\xi\circ f(S(x_{n-1}))
+\frac{1}{2^n}.$$ Then
$$0\leq \inf\,\xi\circ f(x_n) -\inf\,\xi\circ f(S(x_{n-1})) <
\epsilon_n\ \rightarrow\ 0\ \  (n\rightarrow\infty).$$ For each $n$,
take $y_n\in f(x_n)$ such that
$$\xi\circ y_n -\inf\,\xi\circ f(S(x_{n-1}))\, <\, \epsilon_n.\eqno{(3.22)}$$
Since $x_n\in S(x_{n-1}),\ \forall n,$ it is easy to see that
$$S(x_0)\supset S(x_1)\supset S(x_2)\supset \cdots.$$
When $m>n$, $x_m\in S(x_{m-1})\subset S(x_n)$. Thus,
$$y_n\in f(x_n) \subset f(x_m) +\frac{\epsilon}{\lambda} d(x_m, x_n)
k_0 +D.$$ Hence, there exists $y^{\prime}_m\in f(x_m)$ such that
$$y_n -y_0\,\geq_D\, y^{\prime}_m -y_0 +\frac{\epsilon}{\lambda}
d(x_m, x_n) k_0.\eqno{(3.23)}$$ From (3.23) and using Lemma 3.2, we
have
$$\xi_{k_0}(y_n-y_0)\,\geq\, \xi_{k_0}(y^{\prime}_m -y_0)
+\frac{\epsilon}{\lambda} d(x_m, x_n).$$ Remarking that
$y^{\prime}_m\in f(x_m)\subset f(S(x_{n-1}))$ and using (3.22),  we
have
\begin{eqnarray*}
\frac{\epsilon}{\lambda} \,d(x_m, x_n) &\leq& \xi_{k_0}(y_n-y_0) -
\xi_{k_0}(y^{\prime}_m -y_0)\\
&=& \xi(y_n) -\xi(y^{\prime}_m)\\
&\leq& \xi\circ y_n -\inf\,\xi\circ f(S(x_{n-1}))\\
&<& \epsilon_n\ \rightarrow \ 0\ \ (m>n\rightarrow\infty).
\end{eqnarray*}
From this, we conclude that $(x_n)$ is a Cauchy sequence  which
satisfies that $S(x_{n+1})\subset S(x_n) \subset\cdots\subset
S(x_0),\ \forall n.$ As $(X, d)$ is $S(x_0)$-dynamically  complete,
there  exists $u\in X$ such that $x_n\rightarrow u$. Let $n$ be
given. Then for any $i\in N$, $x_{n+i}\in S(x_n)$ and
$$f(x_n) \subset f(x_{n+i}) +\frac{\epsilon}{\lambda} d(x_{n+i},
x_n) k_0 +D. \eqno{(3.24)}$$ Since $f$ is $D$-s.l.m. and $f$ has
$D$-$k_0$-closed values, i.e.,  $f(x)+D$ is $k_0$-closed for any
$x\in X$, from (3.24) we can deduce that
$$f(x_n) \subset f(u) +\frac{\epsilon}{\lambda} d(u, x_n) k_0 +D.$$
That is, $u\in S(x_n)$ and assumption (F) is satisfied.  Thus, we
can apply Theorem 3.1 and obtain $\hat{x}\in X$ such that (a) and
(b) hold. As done in the proof of Corollary 3.5, we can also deduce
that (c) holds.
\hfill\framebox[2mm]{}\\

{\bf Remark 3.1.} \   In order to compare Corollary 3.6 with [41,
Theorem 3.1], let us recall the notion of $(f, D)$-lower
completeness. A metric space $(X, d)$ is said to be $(f, D)$-lower
complete if every Cauchy sequence $(x_n)\subset X$ satisfying
$f(x_n)\subset f(x_{n+1}) +D$ for each $n$, is convergent. If a
Cauchy sequence $(x_n)\subset  S(x_0)$ satisfies that $x_{n+1}\in
S(x_n)$, i.e., $ f(x_n) \subset f(x_{n+1}) +(\epsilon/\lambda)
d(x_{n+1}, x_n)  k_0 +D$, then $f(x_n) \subset f(x_{n+1}) +D$.
Hence, it is obvious that $(X, d)$ being $(f, D)$- lower complete
implies that $(X, d)$ is $S(x_0)$-dynamically complete. Besides, in
Corollary 3.6 we only require that $k_0\in D\backslash -{\rm
vcl}(D)$ and $f$ has $D$-$k_0$-closed values,  which are
respectively weaker than the conditions that $k_0\in
D\backslash-{\rm cl}(D)$ and that $f$ has $D$-closed values. Hence,
Corollary 3.6 is indeed a generalization of [41, Theorem 3.1].

Similarly, in Corollaries 3.4 and 3.5, the assumption that $(X, d)$
is complete can also be replaced by a weaker assumption: $(X, d)$ is
$S(x_0)$-dynamically complete.\\

{\section*{ \large\bf 4.   Set-valued EVPs with perturbations
containing a set   }

\hspace*{\parindent} In this section, by using the results in
Section 3 we give several set-valued EVPs, where perturbations
containing a convex subset of the ordering cone. First we give a
generalization of [33,  Theorem 3.4].\\

{\bf Theorem 4.1.}   \ {\sl Let $(X, d)$ be a complete metric space.
$Y$ be a locally convex space pre-ordered by a convex cone $D$,
$H\subset D\backslash\{0\}$ be a convex set and  $f:\, X\rightarrow
2^Y\backslash \{\emptyset\}$ be a set-valued map.

Suppose that the following assumptions are satisfied:

{\rm (B$_1$)} \ $0\not\in {\rm cl}(H+D)$.

{\rm (B$_2$)} \ $f(X)$ is $D$-bounded.

{\rm (B$_3$)} \ $f$ is $D$-s.l.m.

Then for any $x_0\in X$ and any $\gamma >0$, there exists
$\hat{x}\in X$ such that

{\rm (a)}  $f(x_0)\subset f(\hat{x}) +\gamma^{\prime} d(\hat{x},
x_0)H +D,\ \ \forall\gamma^{\prime}\in (0,\gamma)$;

{\rm (b)} $\forall x\in X\backslash\{\hat{x}\},\ \exists
\gamma^{\prime}\in (0,\gamma)$ such that $f(\hat{x})\not\subset f(x)
+\gamma^{\prime} d(x, \hat{x}) H +D.$}\\

{\bf Proof.}\  Put $F_{\gamma^{\prime}}(x,
x^{\prime}):=\gamma^{\prime} d(x, x^{\prime}) H,\
\forall\gamma^{\prime}\in (0,\gamma)$. Clearly, the family
$\{F_{\gamma^{\prime}}\}_{\gamma^{\prime}\in (0,\gamma)}$ satisfies
the property TI. Put $$S(x_0):=\{x\in X:\, f(x_0)\subset f(x)
+\gamma^{\prime} d(x, x_0) +D,\ \forall \gamma^{\prime}\in (0,
\gamma)\}.$$ Obviously, $x_0\in S(x_0)$ and $S(x_0)\not=\emptyset$.
By assumption (B$_1$) and the separation theorem, there exists
$\xi\in Y^*$ and $\alpha >0$ such that $\xi(H+D)\geq \alpha>0$.
Thus, $\xi\in D^+\backslash\{0\}$ and $\xi(H)\geq\alpha>0$. By
assumption (B$_2$), $\xi$ is lower bounded on $f(X)$ and lower
bounded on $f(S(x_0))$. Take any fixed $\gamma_0\in (0,\gamma)$. For
any $\delta>0$, put $F_{\gamma_0\delta}:=\bigcup\{\gamma_0 d(x,
x^{\prime})H:\, d(x, x^{\prime})\geq\delta\}$. Obviously,
$$\inf\xi\circ F_{\gamma_0\delta}\,\geq\,\gamma_0\,\delta\,
\alpha\,>\,0.$$ Hence, assumption (E$_3$) in Corollary 3.4 is
satisfied. In order to apply Corollary 3.4 we need to show that for
any $x\in S(x_0)$, $S(x)$ is dynamically closed. Let $(x_n)\subset
S(x)$ such that $x_{n+1}\in S(x_n)\subset S(x), \ \forall n$, and \
$x_n\rightarrow u$. For any $n\in N$ and any $\gamma^{\prime}\in
(0,\gamma)$, we have
$$ f(x_n)\subset f(x_{n+1}) +\gamma^{\prime} d(x_{n+1}, x_n)
H+D\subset f(x_{n+1}) +D.$$ By assumption (B$_3$) we have
$$f(x_n)\,\subset\, f(u) +D.\eqno{(4.1)}$$
For any $k\geq n$, $x_k\in S(x_n)$.  By this and (4.1), for any
$\gamma^{\prime\prime}\in (0,\gamma)$, we have
$$f(x_n)\,\subset\, f(x_k) +{\gamma}^{\prime\prime} d(x_k, x_n)
H+D\,\subset\, f(u) +\gamma^{\prime\prime} d(x_k,
x_n)H+D.\eqno{(4.2)}$$ For any fixed $\gamma^{\prime}\in (0,
\gamma)$, take any $\gamma^{\prime\prime}\in
(\gamma^{\prime},\gamma)$. Since $d(x_k, x_n)\rightarrow d(u, x_n)\
(k\rightarrow\infty)$, there exists $k^{\prime}\geq n$ such that
$$d(x_{k^{\prime}}, x_n)\,\geq\,
\frac{\gamma^{\prime}}{\gamma^{\prime\prime}}\, d(u, x_n),\ \ {\rm
i.e.,}\ \ \gamma^{\prime\prime}\, d(x_{k^{\prime}}, x_n)\,\geq\,
\gamma^{\prime}\, d(u, x_n). \eqno{(4.3)}  $$ By (4.2) and (4.3), we
have
\begin{eqnarray*}
f(x_n)\, &\subset & \, f(u) +\gamma^{\prime\prime} d(x_{k^{\prime}},
x_n) H+D\\
&\subset& \, f(u) +\gamma^{\prime} d(u, x_n) H+D.
\end{eqnarray*}
Thus, $u\in S(x_n)\subset S(x)$ and $S(x)$ is dynamically closed.
Now, applying Corollary 3.4 we obtain the result. \hfill\framebox[2mm]{}\\

{\bf Remark 4.1.} \ If we denote the set $\{l\in Y^*:\,\inf l\circ
H>0\}$ by $H^{+s}$, then (B$_1$) is equivalent to the following
assumption:

(B$_1^{\prime}$) \ $H^{+s}\cap D^+\not=\emptyset$.

In this expression, (B$_2$) can be relaxed to the following weaker
assumption:

(B$_2^{\prime}$) \ $\exists\, \xi\in H^{+s}\cap D^+$ such that $\xi$
is lower bounded on $f(X)$.\\

{\bf Remark 4.2.} \  When $Y$ is a Banach space, [33, Theorem 3.4]
gave the same result under the following assumptions (A1)-(A4):

(A1) \ $H\subset D\backslash\{0\}$ is a closed convex set (thus we
have $\kappa:=d(0, H) >0$).

(A2) \ $\zeta:=\inf \{d(h/\|h\|, -D):\, h\in H\}\,>\,0$.

(A3) \  $f(X)$ is $D$-bounded.

(A4) \ ${\rm epi} f$ is closed in $X\times Y$, where ${\rm epi}
f\,=\,\{(x, y)\in X\times Y:\,y\in f(x)+D\}$.

Obviously, assumption (A3) is exactly (B$_2$) here. It is easy to
show that assumptions (A1) and (A2) imply (B$_1$). Here, we needn't
assume that $H$ is closed. The essential assumption is $0\not\in
{\rm cl}(H+D)$. Next we show that assumption (B$_3$) is strictly
weaker than (A4). It is easy to see that (A4) $\Rightarrow$ (B$_3$).
In fact, assume that (A4) holds. Let $(x_n)\subset X$ such that
$f(x_n)\subset f(x_{n+1}) +D,\ \forall n,$ and $x_n\rightarrow u$.
Let $n$ be given. For any $k\geq n$, we have $f(x_n)\subset f(x_k)
+D$. Take any $y_n\in f(x_n)$. Then the points $(x_k, y_n)\in {\rm
epi} f$. Since $(x_k, y_n)\rightarrow (u, y_n)\
(k\rightarrow\infty)$ in $X\times Y$ and ${\rm epi} f$ is closed, we
have $y_n\in f(u) +D$ and hence $f(x_n)\subset f(u) +D$. That is,
$f$ is $D$-s.l.m. The following example shows that there exists a
$D$-s.l.m. set-valued map $f$ such that ${\rm epi} f$ is not closed,
i.e., (B$_3$)$\not\Rightarrow$ (A4). Hence, even in the case that
$Y$ is a Banach space, Theorem 4.1 is also an improvement of [33,
Theorem 3.4].\\
\\

{\bf Example 4.1.} \  A $D$-s.l.m. set-valued map $f$ such that
${\rm epi} f$ is not closed. Let $X=R$, $Y=R$ and $D=[0,+\infty)$.
Define a set-valued map $f:\,  X\rightarrow
2^Y\backslash\{\emptyset\}$ as follows:
$$
f(x)=\left\{
\begin{array}{cc}
x+1+D, \ \ &\ \ {\rm if}\ \ x\in [0,+\infty),\\
\{x\}, \ \ &\ \ {\rm if}\ \ x\in(-\infty, 0).
\end{array}
\right.
$$
It is easy to verify that $f$ is $D$-s.l.m. However, ${\rm epi} f$
is not closed in $X\times Y$. For example, take a sequence $(x_n)$
in the interval $(-\infty, 0)$ such that $x_n\rightarrow  0$.
Obviously, $(x_n, x_n)\rightarrow (0, 0)$ in $X\times Y= R^2$. Here,
every $(x_n, x_n)\in {\rm epi} f$. But $0\not\in f(0) +D$, i.e.,
$(0, 0)\not\in
{\rm epi} f$. Thus, ${\rm epi} f$ is not closed.\\

Next we further introduce the following assumption:

(B$_3^{\prime}$) \ $f$ is $D$-s.l.m. and has $D$-closed values.

It is easy to see that (A4) $\Rightarrow$ (B$_3^{\prime}$)
$\Rightarrow$ (B$_3$). In Example 4.1, $f$ also has $D$-closed
values. Hence Example 4.1 shows that (B$_3^{\prime}$)
$\not\Rightarrow$ (A4). The following theorem  is a set-valued
extension of both [48, Theorem 6.2] and [42, Theorem 6.8], and it
also generalizes and improves [33, Theorem 3.5] in the case that
(ii) holds, where (ii) means that $H$ is bounded. First let us
recall some concepts (see [37, Definition 2.1.4]). Let $H$ be a
subset of a topological vector space. A convex series of points in
$H$ is a series of the form $\sum_{n=1}^{\infty}\lambda_n x_n$,
where $x_n\in H,\,\lambda_n\geq 0$ and $\sum_{n=1}^{\infty}\lambda_n
=1$. $H$ is said to be $\sigma$-convex if every convex series of its
points converges to a point of $H$. Sometimes, a $\sigma$-convex set
is called a cs-complete set, see, e.g. [48, 51]. In [33, Theorem
3.5], $H\subset D\backslash\{0\}$ is assumed to be a closed convex
subset of a Banach space $Y$ and (ii) means that $H$ is bounded, so
$H$ there is a $\sigma$-convex set. But a $\sigma$-convex set may be
non-closed. For example, an open ball in a Banach space is a
$\sigma$-convex set but non-closed. For details, see e.g., [43] and
the references therein. As every singleton is $\sigma$-convex, if we
take $H=\{k_0\}$,
where $k_0\in D\backslash -{\rm cl}(D)$, then we can also deduce Corollary 3.5 from the following theorem.\\

{\bf Theorem 4.2.}   \ {\sl Let $(X, d)$ be a complete metric space,
$x_0\in X$,  $Y$ be a locally convex space quasi ordered by a convex
cone $D$, $H\subset D\backslash\{0\}$ be a convex set and $f:\,
X\rightarrow 2^Y\backslash\{\emptyset\}$ be a set-valued map.
Suppose that the following assumptions are satisfied:

{\rm (B$_1^{\prime}$)} \ $H^{+s}\cap D^+\not= \emptyset$, or
equivalently, $0\not\in {\rm cl}(H+D)$.

{\rm (B$_2^{\prime}$)} \ $\exists\,\xi\in H^{+s}\cap D^+$ such that
$\xi$ is lower bounded on $f(S(x_0))$, where $S(x_0)=\{x\in X:\,
f(x_0)\subset f(x) +\gamma d(x, x_0)H +D\}$.

{\rm (B$_3^{\prime}$)} \ $f$ is $D$-s.l.m. and has $D$-closed
values.

Moreover, suppose that $H$ is $\sigma$-convex, or, $Y$ is locally
complete and $H$ is  locally closed, bounded.

Then for  any $\gamma >0$, there exists $\hat{x}\in X$ such that

{\rm (a)}  $f(x_0)\subset f(\hat{x}) +\gamma d(\hat{x}, x_0)H +D;$

{\rm (b)} $\forall x\in X\backslash\{\hat{x}\}, \
f(\hat{x})\not\subset f(x) +\gamma  d(x, \hat{x}) H +D.$

 {\rm Concerning local completeness and local closedness, see [37,
Chapter 5] and [38, 39, 46]}}
\\

{\bf Proof.}\ Put $F(x, x^{\prime}) = \gamma d(x, x^{\prime}) H$.
Obviously, the family $\{F\}$ satisfies the property TI. Also, it is
obvious that $x_0\in S(x_0)$ and $S(x_0)\not=\emptyset$. By
(B$_1^{\prime}$) and (B$_2^{\prime}$), we can easily show that
assumptions (D) and (E$_3$) in Corollary 3.4 are satisfied. In order
to apply Corollary 3.4, we only need to show that for any $x\in
S(x_0)$, $S(x)$ is dynamically closed. Let $(x_n)\subset S(x),\
x_{n+1}\in S(x_n)\subset S(x),\ \forall n,$ and $x_n\rightarrow u$.
Take any fixed $n_0\in N$ and put $z_1:=x_{n_0}$. As $d(x_k,
u)\rightarrow 0\ (k\rightarrow \infty)$, we may choose a sequence
$(z_n)$ from $(x_k)$ such that $d(z_{n+1}, u) < 1/(n+1)$ and
$z_{n+1}\in S(z_n),\ \forall n$. Take any $y_1\in f(z_1)$. As
$z_2\in S(z_1)$, we have
$$y_1\in f(z_1)\subset f(z_2) +\gamma d(z_2, z_1)H +D.$$
Hence, there exists $y_2\in f(z_2),\ h_1\in H$ and $d_1\in D$ such
that
$$y_1=y_2 +\gamma d(z_2, z_1) h_1 +d_1.$$
In general, if $y_n\in f(z_n)$ is given, then
$$y_n\in f(z_n)\subset f(z_{n+1})+\gamma d(z_{n+1}, z_n) H +D,$$
so there exists $y_{n+1}\in f(z_{n+1}),\ h_n\in H$ and $d_n\in D$
such that
$$y_n = y_{n+1} +\gamma d(z_{n+1}, z_n) h_n +d_n.$$
Adding two sides of the above $n$ equalities, we have
$$\sum\limits_{i=1}^n y_i\,=\,\sum\limits_{i=2}^{n+1} y_i +\gamma
\sum\limits_{i=1}^n d(z_{i+1}, z_i) h_i +\sum\limits_{i=1}^n d_i.$$
From this,
$$y_1\,=\, y_{n+1} +\gamma\sum\limits_{i=1}^n d(z_{i+1}, z_i) h_i
+\sum\limits_{i=1}^n d_i.\eqno{(4.4)}$$

As $\xi\in H^{+s}\cap D^+$, $\xi(D)\geq 0$ and there exists $\alpha
>0$ such that $\xi(H)\geq\alpha$. Acting on two sides of (4.4) by
$\xi$, we have
\begin{eqnarray*}
\xi\circ y_1\,&=&\, \xi\circ y_{n+1} +\gamma \sum\limits_{i=1}^n
d(z_{i+1}, z_i)\,\xi(h_i)  + \xi\left(\sum\limits_{i=1}^n
d_i\right).\\
&\geq&\, \xi\circ y_{n+1} +\gamma\alpha\left(\sum\limits_{i=1}^n
d(z_{i+1}, z_i)\right).
\end{eqnarray*}
From this and by (B$_2^{\prime}$),
\begin{eqnarray*}
\sum\limits_{i=1}^n d(z_{i+1},
z_i)\,&\leq&\,\frac{1}{\gamma\alpha}(\xi\circ y_1-\xi\circ
y_{n+1})\\
&\leq&\,\frac{1}{\gamma\alpha}(\xi\circ y_1-\inf \xi\circ
f(S(x_0)))\\
&<&\,+\infty.
\end{eqnarray*}
Hence, $\sum\limits_{i=1}^{\infty}d(z_{i+1}, z_i)\,<\, +\infty$. By
the assumption that $H$ is $\sigma$-convex, we conclude that
$$\frac{\sum\limits_{i=1}^{\infty} d(z_{i+1}, z_i)
h_i}{\sum\limits_{j=1}^{\infty} d(z_{j+1}, z_j)} $$  is  convergent
to  some  point    $\bar{h}\in H.$ Put
$$h_n^{\prime}:=\frac{\sum\limits_{i=1}^n d(z_{i+1}, z_i)
h_i}{\sum\limits_{j=1}^n d(z_{j+1}, z_j)}.$$ Then every
$h_n^{\prime}\in H$ and $h_n^{\prime}\rightarrow\bar{h}$. From
(4.4), we have
$$y_1\,\in\,y_{n+1} +\gamma\left(\sum\limits_{i=1}^n d(z_{i+1},
z_i)\right) h_n^{\prime} +D.\eqno{(4.5)}$$ Remark that
$$\sum\limits_{i=1}^n d(z_{i+1}, z_i)\,\geq\, d(z_1, u) -d(z_{n+1},
u) \ \ {\rm  and}\ \  d(z_{n+1}, u) <1/(n+1).\eqno{(4.6)}$$ Also, By
(B$_3^{\prime}$), $$y_{n+1}\,\in\, f(z_{n+1})\,\subset\,  f(u)
+D.\eqno{(4.7)}$$ Combining (4.5), (4.6) and (4.7), we have
\begin{eqnarray*}
y_1\,&\in&\, y_{n+1} +\gamma(d(z_1, u)-d(z_{n+1}, u)) h_n^{\prime}
+D\\
&\subset&\, y_{n+1} +\gamma (d(z_1, u) -1/(n+1))h_n^{\prime}
+D\\
&\subset&\, f(u) +\gamma (d(z_1, u) - 1/(n+1)) h_n^{\prime}
+D.\hspace{5.6cm} (4.8)
\end{eqnarray*}
Since $\gamma(d(z_1, u)-1/(n+1))h_n^{\prime}\,\rightarrow\, \gamma
d(z_1, u) \bar{h}$ and $f(u) +D$ is closed by (B$_3^{\prime}$), from
(4.8) we have
$$y_1\,\in\, f(u) +\gamma d(z_1, u) \bar{h} +D,\ \ {\rm where}\ \
\bar{h}\in H.$$ Thus, we have shown that
$$f(z_1)\,\subset\, f(u) +\gamma d(z_1, u) H+ D,\ \ {\rm that\ is,}\
\ u\in S(z_1)\subset S(x).$$
 Now, we can apply Corollary 3.4 and the result follows.
Finally, we point out that $Y$ is locally complete iff it is
$l^1$-complete (see [46]). Hence it is easy to see that $Y$ being
locally complete and $H$ being locally closed bounded imply that $H$
is $\sigma$-convex. \hfill\framebox[2mm]{}\\

In Theorem 4.2, if we strengthen assumption (B$_3^{\prime}$) to the
following (B$_3^{\prime\prime}$): $f$ is $D$-s.l.m. and has $(H,
D)$-closed values, i.e., for any $x\in X$ and any $\lambda\geq 0$,
$f(x)+\lambda\,H+D$ is closed, then the assumption that $H$ is
$\sigma$-convex can be weakened to that $H$ is bounded.\\

{\bf Theorem 4.2$^{\prime}$.}   \ {\sl Let $(X, d)$ be a complete
metric space, $x_0\in X$, $Y$ be a locally convex space pre-ordered
by a convex cone $D$, $H\subset D\backslash\{0\}$ be a convex set
and $f:\, X\rightarrow 2^Y\backslash\{\emptyset\}$ be a set-valued
map. Suppose that the following assumptions are satisfied:

{\rm (B$_1^{\prime}$)} \ $H^{+s}\cap D^+\not= \emptyset$, or
equivalently, $0\not\in {\rm cl}(H+D)$.

{\rm (B$_2^{\prime}$)} \ $\exists\,\xi\in H^{+s}\cap D^+$ such that
$\xi$ is lower bounded on $f(S(x_0))$.

{\rm (B$_3^{\prime\prime}$)} \ $f$ is $D$-s.l.m. and has $(H,
D)$-closed values.

Moreover, suppose that $H$ is bounded.

 Then for any $\gamma
>0$, there exists $\hat{x}\in X$ such that

{\rm (a)}  $f(x_0)\subset f(\hat{x}) +\gamma d(\hat{x}, x_0)H +D;$

{\rm (b)} $\forall x\in X\backslash\{\hat{x}\}, \
f(\hat{x})\not\subset f(x) +\gamma  d(x, \hat{x}) H +D.$ }
\\

{\bf Proof.}\  As shown in the proof of Theorem 4.2, we easily see
that assumption (E$_3$) in Corollary 3.4 is satisfied. In order to
apply Corollary 3.4, we only need to show that for any $x\in
S(x_0)$, $S(x)$ is dynamically closed. Let $(x_n)\subset S(x)$,
$x_{n+1}\in S(x_n)\subset S(x),\ \forall n$, and $x_n\rightarrow u$.
Take any fixed $n_0\in N$ and put $z_1:= x_{n_0}$. As $d(x_k,
u)\rightarrow 0\ (k\rightarrow\infty)$, we may choose a sequence
$(z_n)$ from $(x_k)$ such that $d(z_{n+1}, u) < 1/(n+1)$ and
$z_{n+1}\in S(z_n),\ \forall n$. Take any $y_1\in f(z_1)$. As done
in the proof of Theorem 4.2, for every $n$, we can choose
$y_{n+1}\in f(z_{n+1}),\, h_n\in H$ and $d_n\in D$ such that
$$y_1\,=\, y_{n+1} +\gamma \left(\sum\limits_{i=1}^n d(z_{i+1},
z_i)h_i\right) +\sum\limits_{i=1}^n d_i\,=\, y_{n+1} +\gamma
\left(\sum\limits_{i=1}^n d(z_{i+1}, z_i)\right) h_n^{\prime}
+\sum\limits_{i=1}^n d_i,$$ where
$$h_n^{\prime}\,=\, \frac{\sum\limits_{i=1}^n d(z_{i+1}, z_i)
h_i}{\sum\limits_{j=1}^n d(z_{j+1}, z_j)} \,\in\, H.$$ Combining
this with the assumption that $d(z_{n+1}, u) < 1/(n+1)$, we can also
deduce (4.8), i.e.,
$$y_1\in f(u) +\gamma\left( d(z_1, u) -\frac{1}{n+1}\right)
h_n^{\prime} +D.$$ From this,
$$y_1+\frac{\gamma}{n+1} h_n^{\prime}\,\in\, f(u) +\gamma d(z_1, u)
H+D.$$ Since $f(u) +\gamma d(z_1, u)H+D$ is closed by
(B$_3^{\prime\prime}$), letting $n\rightarrow\infty$ we have
$$y_1\,\in\, f(u) +\gamma d(z_1, u)H +D.$$
Thus,
$$f(z_1)\subset f(u) +\gamma d(z_1, u)H+D\ \ {\rm and}\ \
u\in S(z_1)= S(x_{n_0})\subset S(x).$$ Now, applying Corollary 3.4,
we obtain the result. \hfill\framebox[2mm]{}\\

{\bf Theorem 4.3.}   \ {\sl Let $(X, d)$ be a complete metric space,
$Y$ be a locally convex space whose topology is determined by a
saturated family $\{p_{\alpha}\}_{\alpha\in\Lambda}$ of semi-norms
{\rm (concerning saturated family of semi-norms, see [26, p.96])},
$D\subset Y$ be a convex cone and $H\subset D\backslash\{0\}$ be a
closed convex set. Suppose that the following assumptions are
satisfied:

{\rm (B$_2$)}  \   $f(X)$ is $D$-bounded.

{\rm (B$_3^{\prime}$)} \  $f$ is $D$-s.l.m. and has $D$-closed
values.

Moreover, assume that $Y$ is $l^{\infty}$-complete {\rm (see [37,
38])} and for each $\alpha\in\Lambda$ there exists $\xi_{\alpha}\in
D^+\backslash\{0\}$ and $\lambda_{\alpha}>0$ such that
$\lambda_{\alpha}p_{\alpha}(h)\leq\xi_{\alpha}(h),\ \forall h\in H.$

Then for any $x_0\in X$ and any $\gamma >0$, there exists
$\hat{x}\in X$ such that

{\rm (a)}  $f(x_0)\subset f(\hat{x}) +\gamma d(\hat{x}, x_0)H +D;$

{\rm (b)} $\forall x\in X\backslash\{\hat{x}\}, \
f(\hat{x})\not\subset f(x)
+\gamma  d(x, \hat{x}) H +D.$}\\

{\bf Proof.}\  Put $F(x, x^{\prime})=\gamma d(x, x^{\prime}) H$.
Obviously, the family $\{F\}$ satisfies the property TI. Also,
$x_0\in S(x_0)$ and $S(x_0)\not=\emptyset$, where $S(x_0)$ is the
same as one in the proof of Theorem 4.2. Since $H$ is closed and
$0\not\in H$, there exists $\alpha_0\in\Lambda$ and $\eta >0$ such
that $p_{\alpha_0}(h)\geq\eta,\ \forall h\in H$. By the assumption,
there exists $\xi_{\alpha_0}\in D^+\backslash\{0\}$ and
$\lambda_{\alpha_0}>0$ such that
$$\lambda_{\alpha_0}\eta\,\leq\,\lambda_{\alpha_0}\,p_{\alpha_0}(h)\,\leq\,\xi_{\alpha_0}(h),\
\ \forall h\in H. \eqno{(4.9)}$$ Thus, $\xi_{\alpha_0}\in D^+\cap
H^{+s}$. Combining this with (B$_2$), we can show that assumption
(E$_3$) in Corollary 3.4 is satisfied. By Corollary 3.4, it is
sufficient to show that for any $x\in S(x_0)$, $S(x)$ is dynamically
closed. Let $(x_n)\subset S(x)$, $x_{n+1}\in S(x_n)\subset S(x),\
\forall n, $ and $x_n\rightarrow u$. As done in the proof of Theorem
4.2, we can obtain a sequence $(z_n)$ from $(x_k)$ such that
$d(z_{n+1}, u) < 1/(n+1)$ and $z_{n+1}\in S(z_n),\ \forall n$. Take
any $y_1\in f(z_1)$. We may choose $y_{n+1}\in f(z_{n+1}),\ h_n\in
H$ and $d_n\in D$ such that (see (4.4))
$$y_1\,=\, y_{n+1} +\gamma \sum\limits_{i=1}^n d(z_{i+1}, z_i) h_i +
\sum\limits_{i=1}^n d_i.\eqno{(4.10)}$$ For each $\alpha\in\Lambda$,
acting on two sides of (4.10) by $\xi_{\alpha}$, we have
\begin{eqnarray*}
\xi_{\alpha}\circ y_1\,&=&\, \xi_{\alpha}\circ y_{n+1}
+\gamma\sum\limits_{i=1}^n d(z_{i+1}, z_i)\,\xi_{\alpha}\circ
h_i+\sum\limits_{i=1}^n\xi_{\alpha}\circ d_i\\
&\geq&\,\xi_{\alpha}\circ y_{n+1} +\gamma\sum\limits_{i=1}^n
d(z_{i+1}, z_i)\,\xi_{\alpha}\circ h_i\\
&\geq&\, \xi_{\alpha}\circ y_{n+1} +\gamma
\lambda_{\alpha}\sum\limits_{i=1}^n d(z_{i+1}, z_i) p_{\alpha}(h_i).
\end{eqnarray*}
From this,
\begin{eqnarray*}
\sum\limits_{i=1}^n d(z_{i+1}, z_i)
p_{\alpha}(h_i)\,&\leq&\,\frac{1}{\gamma\lambda_{\alpha}}(\xi_{\alpha}\circ
y_1-\xi_{\alpha}\circ y_{n+1})\\
&\leq&\,\frac{1}{\gamma \lambda_{\alpha}}(\xi_{\alpha}\circ y_1-\inf
\xi_{\alpha}\circ f(S(x_0)))\\
&<&\,+\infty.
\end{eqnarray*}
Hence, $$\sum\limits_{i=1}^{\infty} d(z_{i+1}, z_i)
p_{\alpha}(h_i)\,<\,+\infty,\ \ \forall \alpha\in\Lambda.$$ Since
$Y$ is $l^{\infty}$-complete, we know that
$\sum\limits_{i=1}^{\infty} d(z_{i+1}, z_i) h_i$ is convergent in
$Y$. Combining this with (4.9), we have
$$\lambda_{\alpha_0}\eta\sum\limits_{i=1}^n d(z_{i+1},
z_i)\,\leq\,\sum\limits_{i=1}^n d(z_{i+1},
z_i)\,\xi_{\alpha_0}(h_i)\,\leq\,
\xi_{\alpha_0}\left(\sum\limits_{i=1}^{\infty} d(z_{i+1}, z_i)
h_i\right),\ \ \forall n\in N.$$ Thus, $$\sum\limits_{i=1}^{\infty}
d(z_{i+1}, z_i)\,<\,+\infty.$$  Put
$$h_n^{\prime}:=\,\frac{\sum\limits_{i=1}^n d(z_{i+1}, z_i)
h_i}{\sum\limits_{j=1}^n d(z_{j+1}, z_j)}.$$ Then
$$ h_n^{\prime}\in H\ \ \ {\rm and}\ \ \ h_n^{\prime} \rightarrow
\bar{h}:=\,\frac{\sum\limits_{i=1}^{\infty} d(z_{i+1}, z_i)
h_i}{\sum\limits_{j=1}^{\infty} d(z_{j+1}, z_j)}\,\in\, H.$$ From
(4.10), we have
$$y_1\,\in\, y_{n+1}+\gamma\left(\sum\limits_{i=1}^n d(z_{i+1},
z_i)\right) h_n^{\prime} +D.$$ This is exactly (4.5). The remains of
the proof is the same as one in the proof of Theorem 4.2 and we omit
the details.\hfill\framebox[2mm]{}\\

{\bf Remark 4.3.} \ Here we needn't assume that (B$_1^{\prime}$) (
equivalently, (B$_1$)) holds in advance. It can be deduced from $H$
being closed and the existence of $\xi_{\alpha}\in D^+\backslash
\{0\}$ and $\lambda_{\alpha}>0$ such that $\lambda_{\alpha}
p_{\alpha}(h)\leq \xi_{\alpha}(h),\ \forall h\in H$. In fact, in
[33, Theorem  3.5], assumption (A2) can be removed when(i) holds,
where $Y$ is a Banach space and (i) means that there exists $\xi\in
D^+\backslash \{0\}$ and $\lambda >0$ such that $\lambda \|h\|\leq
\xi(h),\ \forall h\in H.$ Obviously, Theorem 4.3
is a generalization of the part (where (i) holds) of [33, Theorem 3.5].\\

Similarly, by Corollary 3.4 we can also obtain the part (where (iii)
holds, i.e., $Y$ is reflexive) of [33,  Theorem 3.5]. In this case,
assumption (A4) (i.e., ${\rm epi}f$ is closed in $X\times Y$) in
[33] can also be relaxed to (B$_3^{\prime}$) (i.e., $f$ is
$D$-s.l.m. and has $D$-closed values). We already presented several
set-valued EVPs, where perturbations are given by a convex subset
multiplied by the distance. Next we further consider more general
version of set-valued EVP, where the perturbation is given by a
convex subset multiplied by a real function which is more
general than the distance.\\

{\bf Theorem 4.4.}   \ {\sl Let $(X, d)$ be a metric space, $Y$ be a
locally convex space, $D\subset Y$ be a convex cone, $H\subset
D\backslash\{0\}$ be a convex set and $f:\, X\rightarrow
2^Y\backslash\{\emptyset\}$ be a set-valued map. Let a real function
$p: X\times X \rightarrow R^+:=[0,+\infty)$ satisfy the following
properties:

{\rm (p$_1$)} \ for any $x_1,\, x_2,\, x_3\in X,\ p(x_1, x_3)\leq
p(x_1, x_2) + p(x_2, x_3)$;

{\rm (p$_2$)} \ every sequence $(x_n)$ with $p(x_n, x_m)\rightarrow
0\ (m>n\rightarrow\infty )$ is a Cauchy sequence, where $p(x_n,
x_m)\rightarrow 0\ (m>n\rightarrow\infty)$ means that for any
$\epsilon >0$, there exists $n_0\in N$ such that $p(x_n, x_m)
<\epsilon$ for all $m>n \geq n_0$;

{\rm (p$_3$)} \  $p(x, x^{\prime})>0,\ \forall x\not=x^{\prime}$.

Let $x_0\in X$ such that $$S(x_0):=\,\{x\in X:\, f(x_0)\subset f(x)
+p(x_0, x)H+D\}\not=\emptyset$$ and $(X, d)$ be $S(x_0)$-dynamically
complete. Suppose that for any $x\in S(x_0)$, $S(x)$ is dynamically
closed and the following assumptions are satisfied:

{\rm (B$_1^{\prime}$)} \ $H^{+s}\cap D^+\not=\emptyset$, or
equivalently, $0\not\in {\rm cl}(H+D)$.

{\rm (B$_2^{\prime}$)} \  $\exists\,\xi\in H^{+s}\cap D^+$ such that
$\xi$ is lower bounded on $f(S(x_0))$.

Then there exists $\hat{x}\in X$ such that

{\rm (a)}  $f(x_0)\subset f(\hat{x}) + p(x_0, \hat{x})H +D;$

{\rm (b)} $\forall x\in X\backslash\{\hat{x}\}, \
f(\hat{x})\not\subset f(x)
+ p(\hat{x}, x) H +D.$}\\

{\bf Proof.}\ Put $F(x, x^{\prime}):=p(x^{\prime}, x) H,\ \forall x,
x^{\prime} \in X$. Obviously, the family $\{F\}$ satisfies the
property TI. By assumptions (B$_1^{\prime}$) and (B$_2^{\prime}$),
there exists $\xi\in D^+\backslash\{0\}$ and $\alpha>0$ such that
$\xi(H)\geq\alpha>0$ and $\xi$ is lower bounded on $f(S(x_0))$.
Clearly, $\xi$ satisfies assumption (D) in Theorem 3.1$^{\prime}$.
For any $x\in S(x_0)$ and any $x^{\prime}\in S(x)\backslash\{x\}$,
we have
$$f(x)\,\subset\, f(x^{\prime}) +p(x, x^{\prime})H+D.\eqno{(4.11)}$$
By (p$_3$), $p(x, x^{\prime})\alpha >0$, thus we may take $y\in
f(x)$ such that
$$\xi\circ y\,<\, \inf \xi\circ f(x) + \frac{1}{2}p(x,
x^{\prime})\alpha.\eqno{(4.12)}$$ By (4.11), we have
$$y\,\in \, f(x^{\prime}) + p(x, x^{\prime})H +D.$$
Thus, there exists $y^{\prime}\in f(x^{\prime}),\ h^{\prime}\in H$
and $d^{\prime}\in D$ such that $$y=y^{\prime} +p(x,
x^{\prime})h^{\prime} +d^{\prime}.$$ Acting two sides of the above
equality by $\xi$, we have
\begin{eqnarray*}
\xi\circ y&=&\xi\circ y^{\prime} +p(x, x^{\prime})\,\xi\circ
h^{\prime} +\xi\circ d^{\prime}\\
&\geq&  \xi\circ y^{\prime} +p(x, x^{\prime})\,\alpha\\
&\geq& \inf\xi\circ f(x^{\prime}) + p(x, x^{\prime})\,\alpha.
\end{eqnarray*}
Combining this with (4.12), we have
$$\inf\xi\circ f(x)\,>\, \inf\xi\circ f(x^{\prime}) +\frac{1}{2}
p(x, x^{\prime})\,\alpha\,>\, \inf \xi\circ f(x^{\prime}).$$ Thus,
$\xi$ satisfies assumption (E) in Theorem 3.1$^{\prime}$.

Next we show that assumption (F) is satisfied. Let a sequence
$(x_n)\subset S(x_0)$ such that $x_n\in S(x_{n-1})$ and
$$\inf\xi\circ f(x_n)\,<\,\inf\xi\circ f(S(x_{n-1})) +\epsilon_n,\
\forall n,$$  where $\epsilon_n>0$  and $ \epsilon_n\rightarrow 0.$
For each $n$, take $y_n\in f(x_n)$ such that
$$\xi\circ y_n\,<\, \inf\xi\circ
f(S(x_{n-1}))+\epsilon_n.\eqno{(4.13)}$$ When $m>n$, $x_m\in
S(x_n)$. Hence
$$y_n\,\in\, f(x_n)\,\subset\, f(x_m) +p(x_n, x_m) H+D.$$
Thus, there exists $y_{m,n}\in f(x_m)$, $h_{m,n}\in H$ and
$d_{m,n}\in D$ such that
$$y_n\,=\, y_{m,n} + p(x_n, x_m) h_{m,n}+ d_{m,n}.$$
Acting two sides of the above equality by $\xi$, we have
$$\xi\circ y_n\,=\,\xi\circ y_{m, n} +p(x_n, x_m)\,\xi\circ h_{mn}
+\xi\circ d_{mn}\,\geq\, \xi\circ y_{mn} + p(x_n,
x_m)\,\alpha.\eqno{(4.14)}$$ Observe that $y_{m,n}\in f(x_m)\subset
f(S(x_{n-1}))$. From (4.14) and (4.13), we have
\begin{eqnarray*}
p(x_n, x_m)\,&\leq&\,\frac{1}{\alpha} (\xi\circ y_n-\xi\circ
y_{m,n})\\
&\leq&\,\frac{1}{\alpha}(\xi\circ y_n-\inf\xi\circ f(S(x_{n-1})))\\
&<&\,\frac{1}{\alpha} {\epsilon_n}.
\end{eqnarray*}
Hence $p(x_n, x_m)\rightarrow 0\ (m>n\rightarrow\infty)$. By
property (p$_2$),  $(x_n)$ is a Cauchy sequence. As $(X, d)$ is
$S(x_0)$-dynamically complete, there exists $u\in X$ such that
$x_n\rightarrow u$. Since $S(x_n)$ is dynamically closed, we easily
see that $u\in S(x_n),\ \forall n.$ Thus, assumption (F) in Theorem
3.1$^{\prime}$ is
satisfied. By Theorem 3.1$^{\prime}$, we obtain the result. \hfill\framebox[2mm]{}\\

At the end of this section, we consider EVPs for approximately
efficient solutions.  N\'{e}meth [35] gave the concept of
approximately efficient solutions for vector-valued maps. Here, we
extend the concept to set-valued maps.

Let $X$ be a nonempty set, $Y$ be a real linear space, $D\subset Y$
be a convex pointed cone and $f:\, X\rightarrow
2^Y\backslash\{\emptyset\}$ be a set-valued map. We consider the
following vector optimization problem:
$${\rm Min}\{f(x):\, x\in X\}. \eqno{(4.15)}$$
Moreover, let $\epsilon >0$ and $H\subset D\backslash\{0\}$ be a
convex set.\\

{\bf Definition  4.1.}   \   A point $x_0\in X$ is called an
efficient solution of (4.15) if $f(x_0)\not\subset f(X) +
D\backslash\{0\}$, where $f(X)$ denotes the set $\cup_{x\in X}
f(x)$. A point $x_0\in X$ is called an $(\epsilon, H)$-efficient
solution of (4.15) if $f(x_0)\not\subset f(X) +\epsilon\, H +D$, or
equivalently, there exists $y_0\in f(x_0)$ such that
$(y_0-\epsilon\,
H-D)\cap f(X)\,=\,\emptyset$.\\

{\bf Theorem 4.5.}   \ {\sl In Theorem 4.1, moreover assume that
$x_0$ is an $(\epsilon, H)$-efficient solution of (4.15). Then there
exists $\hat{x}\in X$ such that

{\rm (a)}  $f(x_0)\subset f(\hat{x}) +\gamma^{\prime} d(\hat{x},
x_0)H +D,\ \ \forall\gamma^{\prime}\in (0,\gamma)$;

{\rm (b)} $\forall x\in X\backslash\{\hat{x}\},\ \exists
\gamma^{\prime}\in (0,\gamma)$ such that $f(\hat{x})\not\subset f(x)
+\gamma^{\prime} d(x, \hat{x}) H +D;$

{\rm (c)}  $d(\hat{x}, x_0) \leq \epsilon /\gamma$.
}\\

{\bf Proof.} \ By Theorem 4.1, we conclude that (a) and (b) hold.
Hence, we only need to show that (c) holds. If not, assume that
 $d(\hat{x}, x_0) > \epsilon/\gamma$. Then  $\epsilon/d(\hat{x}, x_0) <
\gamma$.  By (a),
\begin{eqnarray*}
f(x_0)&\subset& f(\hat{x}) + (\epsilon/d(\hat{x}, x_0)) d(\hat{x},
x_0) H +D\\
&=& f(\hat{x}) +\epsilon H +D\\
&\subset& f(X) +\epsilon H +D.
\end{eqnarray*}
This contradicts the assumption that $x_0$ is an $(\epsilon,
H)$-efficient solution of (4.15).\hfill\framebox[2mm]{}\\

Similarly we can prove the following:\\

{\bf Theorem 4.6.}   \ {\sl In Theorems 4.2, 4.2$^{\prime}$ and 4.3,
moreover assume that $x_0$ is an $(\epsilon, H)$-efficient solution
of (4.15). Then there exists $\hat{x}\in X$ such that

{\rm (a)}  $f(x_0)\subset f(\hat{x}) +\gamma d(\hat{x}, x_0)H +D;$

{\rm (b)} $\forall x\in X\backslash\{\hat{x}\}, \
f(\hat{x})\not\subset f(x) +\gamma  d(x, \hat{x}) H +D;$

{\rm (c)} $d(\hat{x}, x_0) < \epsilon /\gamma$.
}\\

\section*{ \large\bf 5.  Pre-orders and minimal points in product spaces }

\hspace*{\parindent}  In this section, by using Theorem 2.1 we
discuss pre-orders and minimal points in product spaces.
Particularly, we obtain several versions of EVP for Pareto
minimizers, which generalize the corresponding results in [4, 5, 31,
48].\\

First we recall some concepts on Pareto minimum. Let $Y$ be a real
linear space with a quasi order $\leq_D$ defined by a convex cone
$D$ in $Y$. Let $B\subset Y$ be nonempty. A point $\bar{y}\in B$ is
called a Pareto minimum of $B$ if $y\in B$ and $y\leq_D \bar{y}$
implies that $\bar{y}\leq_D y$. And $\bar{y}\in B$ is called a
strict Pareto minimum of $B$ if $y\not\leq_D \bar{y}$ for all $y\in
B\backslash\{\bar{y}\}$, i.e., $(B-\bar{y})\cap (-D) = \{0\}$. We
denote by ${\rm Min}^D B$ (resp., ${\rm SMin}^D B$) the sets of all
Pareto minima (resp., strict Pareto minima) with respect to the
order $\leq_D$. In general, we have ${\rm SMin}^D B\subset {\rm
Min}^D B$. If $D$ is pointed, i.e., $D\cap(-D) = \{0\}$, then ${\rm
SMin}^D B = {\rm Min}^D B$. A subset $B$ of $Y$ is said to have the
domination (resp., strict domination) property if, for any $y\in B$,
there exists $\bar{y}\in {\rm Min}^D B$ (resp., $\bar{y}\in {\rm
SMin}^D B$) such that $\bar{y}\leq_D y$.

Moreover, let $(X, d)$ be a metric space and let $f:\, X\rightarrow
2^Y\backslash\{\emptyset\}$ be a set-valued map. A point $\bar{x}\in
X$ is called a Pareto minimizer (resp., strict Pareto minmizer) of
$f$ if there exists $\bar{y}\in f(\bar{x})$ such that $\bar{y}\in
{\rm Min}^D f(X)$ (resp., $\bar{y}\in {\rm SMin}^D f(X)$).

Let $F:\, X\times X \rightarrow 2^D\backslash\{\emptyset\}$ satisfy
the following conditions (see [48]):

{\rm (F1)} \  $0\in F(x, x),\ \forall x\in X$;

{\rm (F2)} \  $F(x_1, x_2) +F(x_2, x_3) \subset F(x_1, x_3) +D, \
\forall x_1,\, x_2,\, x_3\in X$;

{\rm (F3)} \  there exists a $D$-monotone extended real function
$\xi: Y \rightarrow (-\infty, +\infty]$  such that for any $y\in Y$
and any $z\in F(X\times X)$, $\xi(y+z) = \xi(y) +\xi(z)$ and such
that for any $\delta >0$,
$$\zeta(\delta):= \inf \{\xi(y):\, y\in \cup_{d(x,
x^{\prime})\geq\delta} F(x, x^{\prime})\} > 0.$$ Obviously, $\xi(0)
=\xi(0+0)=\xi(0) +\xi(0) = 2\xi(0)$, so $\xi(0)=0$.

As in [48], define a quasi-order $\preceq_F$ on $X\times Y$ as
follows:
$$(x_2, y_2) \preceq_F (x_1, y_1)\ \Longleftrightarrow\ y_1\in y_2
+F(x_2, x_1) +D.$$ Consider a nonempty set $A\subset X\times Y$ and
a point $(x_0, y_0)\in A$. Put
$$ S_F(x_0, y_0) := \{(x, y)\in A:\, (x, y)\preceq_F (x_0, y_0)\}
.$$
Moreover, define a partial order $\preceq_{F^*}$ on $X\times Y$
as follows:
$$
(x_2, y_2)\preceq_{F^*} (x_1, y_1) \ \Longleftrightarrow\ \left\{
\begin{array}{cc}
(x_2, y_2) = (x_1, y_1)\ \ {\rm or}\\
(x_2, y_2)\preceq_F(x_1, y_1)\ {\rm and}\ \xi(y_2-y_0) <
\xi(y_1-y_0).
\end{array}
\right.
$$
Put
$$S_{F^*}(x_0, y_0):=\{(x, y)\in A:\, (x,y)\preceq_{F^*} (x_0,
y_0)\}.$$
 If $\xi$ is strict $D$-monotone, i.e., $y_2\leq_D y_1$ and
$y_2\not= y_1$ imply that $\xi(y_2) < \xi(y_1)$, then the orders
$\preceq_F$ and $\preceq_{F^*}$ are coincident.

By using Theorem 2.1, we can deduce the following minimal point
theorem in product spaces. \\

{\bf Theorem 5.1.}   \ {\sl Let $(X, d)$ be a metric space, $Y$ be a
real linear space, $D\subset Y$ be a convex cone and a set-valued
map $F:\, X\times X \rightarrow 2^D\backslash\{\emptyset\}$ satisfy
(F1)-(F3). Let $A\subset X\times Y$ be a nonempty set and $(x_0,
y_0)\in A$ be given. Suppose that the following conditions are
satisfied:

{\rm (i)} \ for any $\preceq_F$-decreasing sequence $\{(x_n, y_n)\}$
in $S_F(x_0, y_0)$, if $\{x_n\}$ is a Cauchy sequence, then
$\{x_n\}$ is convergent in $X$;

{\rm (ii)} \ $\xi$ {\rm (from in (F3))} is lower bounded on
$P_Y(S_F(x_0, y_0))-y_0$, where $P_Y$ is the projection from
$X\times Y$ on $Y$;

{\rm (iii)} \   for every $\preceq_F$-decreasing sequence $\{(x_n,
y_n)\}$ in $S_F(x_0, y_0)$, if $\{x_n\}$ converges to $x$, then
there exists $y\in Y$ such that $(x, y)\in A$ and $(x, y)\preceq_F
(x_n, y_n),\ \forall n$ {\rm (in [48], (iii) is called (H1))}.

Then, there exists $(\hat{x}, \hat{y})\in A$ such that

{\rm (a)} \ $(\hat{x}, \hat{y}) \preceq_{F^*} (x_0, y_0);$

{\rm (b)} \  $(x, y) \preceq_{F^*} (\hat{x}, \hat{y})\
\Longrightarrow \  (x, y) = (\hat{x}, \hat{y}).$

From this, we conclude that $(\hat{x}, \hat{y})\in A$ satisfies the
following

{\rm (a$^{\prime}$)} \  $y_0\in \hat{y} + F(\hat{x}, x_0) +D;$

{\rm (b$^{\prime}$)}  \  for any $(x, y)\in A$ with $x\not=\hat{x}$.
$\hat{y}\not\in y +F(x, \hat{x}) +D.$}\\

{\bf Proof.}\  Obviously, $(A, \preceq_{F^*})$ is a partial order
set and $S_{F^*}(x_0, y_0)\not=\emptyset$. Define  $\eta:\,(A,
\preceq_{F^*}) \rightarrow (-\infty, +\infty]$ as follows:
$\eta\circ(x, y) = \xi(y-y_0), \ \forall (x, y)\in A$. Clearly,
$\eta$ is monotone with respect to $\preceq_{F^*}$.  By (ii), $\xi$
is lower bounded on $P_Y(S_F(x_0, y_0))-y_0$, Also, $0\in
P_Y(S_F(x_0, y_0))-y_0$ and $\xi(0)=0$. Thus,
$$-\infty \,<\, \inf\{\eta\circ (x, y)=\xi(y-y_0):\, (x, y)\in S_{F^*}(x_0,
y_0)\}\,<\, +\infty.$$ This means that assumption (A) is satisfied
if we regard $(A, \preceq_{F^*})$ as a pre-order set $(X, \preceq)$
in Theorem 2.1. For any $(x, y)\in S_{F^*}(x_0, y_0)$ and any
$(x^{\prime}, y^{\prime})\in S_{F^*}(x, y)\backslash\{(x, y)\}$, we
have $y\in y^{\prime} +F(x^{\prime}, x) +D$ and $\xi(y^{\prime}-y_0)
< \xi(y-y_0)$, that is, $\eta\circ(x^{\prime}, y^{\prime})
\,<\,\eta\circ (x, y)$. Thus, assumption (B) in Theorem 2.1 is
satisfied.

Let a sequence $\{(x_n, y_n)\}\subset S_{F^*}(x_0, y_0)$ satisfy
$(x_n, y_n)\preceq_{F^*} (x_{n-1}, y_{n-1})\ \forall n.$ We shall
show that there exists $(\bar{x}, \bar{y})\in A$ such that
$(\bar{x}, \bar{y}) \preceq_{F^*}(x_n, y_n),  \ \forall n$, i.e.,
assumption (C) in Theorem 2.1 satisfied. If there exists a sequence
$n_1 < n_2 < \cdots$ such that $(x_{n_i}, y_{n_i}) = (x_{n_{i+1}},
y_{n_{i+1}}), \ \forall i,$ then we have $(x_k, y_k) = (x_{n_1},
y_{n_1})$ for all $k\geq n_1$ and the result is trivial. Hence,  we
may assume that $(x_n, y_n) \not= (x_{n-1}, y_{n-1}),\ \forall n$.
From the definition of $\preceq_{F^*}$, we have $y_{n-1}\in y_n +
F(x_n, x_{n-1}) +D$ and
$$\xi(y_n-y_0)\,<\, \xi(y_{n-1}-y_0),\ \ \forall n.\eqno{(5.1)}$$
As $\xi$ is lower bounded on $P_Y(S_F(x_0, y_0))-y_0$,
$\{\xi(y_n-y_0)\}$ is a  lower bounded, decreasing real sequence, so
$\{\xi(y_n-y_0)\}$ is convergent. We assume that
$$\gamma:=\lim_{n\rightarrow\infty}\xi(y_n-y_0)\,\in {\bf
R}.\eqno{(5.2)}$$ First, we assert that $\{x_n\}$ is a Cauchy
sequence. If not, there exists $\delta >0$ and a sequence $n_1 < n_2
<\cdots$ such that $d(x_{n_i}, x_{n_{i+1}})\geq \delta$. Since
$(x_{n_{i+1}}, y_{n_{i+1}})\preceq_{F^*} (x_{n_i}, y_{n_i})$, we
have $y_{n_i}\in y_{n_{i+1}}+F(x_{n_{i+1}}, x_{n_i}) +D$. Thus,
$$\xi(y_{n_i}-y_0) -\xi(y_{n_{i+1}}-y_0)\,\geq\, \zeta(\delta)\,>\, 0,\ \
\forall i,$$ which contradicts (5.2). Thus, $\{x_n\}$ is a Cauchy
sequence and by assumption (i), there exists $\bar{x}\in X$ such
that $x_n \rightarrow \bar{x}$. By assumption (iii), there exists
$\bar{y}\in Y$ such that $(\bar{x}, \bar{y})\in A$ and $(\bar{x},
\bar{y})\preceq_F (x_n, y_n),\ \forall n$. From this, we have
$\xi(\bar{y}-y_0)\leq \xi(y_n-y_0)$ and by (5.1), we have
$\xi(\bar{y}-y_0) < \xi(y_n-y_0)$ for all $n$. Thus, $(\bar{x},
\bar{y})\preceq_{F^*} (x_n, y_n),\ \forall n.$  Now, we can apply
Theorem 2.1 and obtain the result, i.e., there exists $(\hat{x},
\hat{y})\in A$ such that (a) and (b) hold. Obviously, $(a)
\Rightarrow  (a^{\prime})$. Next, we show that $(b)\Rightarrow
(b^{\prime})$. Assume that (b$^{\prime}$) is not true. That is,
there exists $(x, y)\in A$ with $x\not= \hat{x}$ such that
$$\hat{y} \in y+ F(x,\hat{x}) +D.\eqno{(5.3)}$$
Thus, there exists $v\in F(x, \hat{x})$ such that $\hat{y}\in y+v+D$
and hence $\hat{y}-y_0\in y-y_0 +v +D$. As $d(x, \hat{x}) >0$, we
have $\xi(v)>0$ and
$$\xi(\hat{y}-y_0) \geq \xi(y-y_0) +\xi(v) >\xi(y-y_0).$$
Combining this with (5.3), we have $(x, y)\preceq_{F^*}
(\hat{x},\hat{y})$ and $x\not=\hat{x}$, which contradicts (b).
\hfill\framebox[2mm]{}\\

{\bf Remark 5.1.}   \  If $\xi$ appeared in (F3) is a positive
linear functional, i.e., $\xi\in D^+\backslash\{0\}$, then the case
becomes easier. This time, we needn't use $y_0$. For example, we may
define $\preceq_{F^*}$ on $X\times Y$ as follows:
$$
(x_2, y_2)\preceq_{F^*} (x_1, y_1) \ \Longleftrightarrow\ \left\{
\begin{array}{cc}
(x_2, y_2) = (x_1, y_1)\ \ {\rm or}\\
(x_2, y_2)\preceq_F(x_1, y_1)\ {\rm and}\ \xi(y_2) < \xi(y_1).
\end{array}
\right.
$$
Moreover, assumption (ii) in Theorem 5.1 can be written as: $\xi$ is
lower bounded on $P_Y(S_F(x_0, y_0))$.

Obviously, Theorem 5.1 generalizes [48, Theorem 2.1].  In fact, it
also includes [31, Theorem 4.2]. In [31, Theorem 4.2], we assume
that $Y$ is a  locally convex space and $H\subset D$ is a convex set
such that $0\not\in {\rm cl}(H+D)$. By the Hahn-Banach separation
theorem, there exists $\xi\in D^+$ such that
$\alpha:=\inf\{\xi(y):\, y\in H\} > 0$. Put $F(x, x^{\prime}) :=
d(x, x^{\prime}) H,\ \forall x, x^{\prime}\in X$. Then for any
$\delta >0$, $\zeta(\delta):=\inf\{\xi(y):\, y\in \cup_{d(x,
x^{\prime})\geq\delta} d(x, x^{\prime}) H\}\geq \alpha \delta >0.$
It is clear that $F$ satisfies (F1)-(F3). Now, applying Theorem 5.1
we can obtain  [31, Theorem 4.2]. In order to obtain  a result on
strict minimal elements in product orders, we need to strengthen
assumption (ii) in Theorem 5.1 (see [31, Theorem 4.3]).\\

{\bf Theorem 5.2.}   \ {\sl  Impose the assumptions of Theorem 5.1
with {\rm (ii)}  being strengthened as

{\rm (ii$^{\prime}$)} \  $\xi$ is lower bounded on $P_Y(S_F(x_0,
y_0))-y_0$ and for all $x\in P_X(S_F(x_0, y_0))$, the set
$\{y^{\prime}\in Y:\, (x, y^{\prime})\in A\}$ has the strict
domination property.

Then, there exists $(\hat{x}, \hat{y})\in A$ such that

{\rm (a)} \  $y_0\in \hat{y} + F(\hat{x}, x_0) +D$ and $\hat{y}\in
{\rm SMin}^D\{y^{\prime}:\, (\hat{x}, y^{\prime})\in A\}$;

{\rm (b)}  \  for any $(x, y)\in A\backslash\{(\hat{x}, \hat{y})\}$,
$\hat{y}\not\in y +F(x, \hat{x}) +D.$}\\

{\bf Proof.}\  By Theorem 5.1, there exists $(\hat{x}, \tilde{y})\in
A$ such that $$y_0\in \tilde{y} + F(\hat{x}, x_0) +D\eqno{(5.4)}$$
and such that $$\tilde{y}\not\in y+F(x, \hat{x}) +D\ {\rm for\ any}\
(x, y)\in A \ {\rm with}\  x\not=\hat{x}.\eqno{(5.5)}$$  By the
imposed strict domination property, there is $\hat{y}\in {\rm
SMin}^D\{y^{\prime}:\, (\hat{x}, y^{\prime})\in A\}$ such that
$\hat{y}\leq_D \tilde{y}$. Next, we show that $(\hat{x}, \hat{y})$
is a desired element. By (5.4) and $\hat{y}\leq_D \tilde{y}$, we
have
 $$y_0\in\hat{y} +D +F(\hat{x}, x_0) +D = \hat{y} +F(\hat{x}, x_0)
+D.$$  Hence, $(\hat{x}, \hat{y})$ satisfies (a).

Let $(x, y)\in A\backslash\{(\hat{x}, \hat{y})\}$. Assume that
$$\hat{y}\in y+F(x, \hat{x}) +D.\eqno{(5.6)}$$
Then $\tilde{y}\in \hat{y} +D \subset y+F(x, \hat{x}) +D$. By (5.5),
we have  $x=\hat{x}$.  Thus, we have $\hat{y},\, y\in
\{y^{\prime}:\, (\hat{x}, y^{\prime})\in A\}$. By (5.6),   $y\leq_D
\hat{y}$. Since $\hat{y}\in{\rm SMin}^D\{y^{\prime}:\, (\hat{x},
y^{\prime})\in A\}$, we  have $y=\hat{y}$. This leads to $(x, y)
=(\hat{x}, \hat{y})$, a contradiction! \hfill\framebox[2mm]{}\\

Let $f:\, X\rightarrow 2^Y\backslash\{\emptyset\}$ and $(x_0,
y_0)\in {\rm gr} f$, where ${\rm gr} f$ denotes the set $\{(x, y)\in
X\times Y:\, x\in X,\, y\in f(x)\}$. For $(x,y),\, (x^{\prime},
y^{\prime})\in {\rm gr} f$, we define $(x^{\prime},
y^{\prime})\preceq_F (x, y)$ iff $y\in y^{\prime} +F(x^{\prime}, x)
+D$. Denote the set  $\{(x, y)\in {\rm gr} f:\, (x, y)\preceq_F
(x_0, y_0)\}$ by $S_F(x_0, y_0)$. By taking  $A={\rm gr} f$ in
Theorem 5.2, we obtain the following.\\

{\bf Corollary 5.1.} \ {\sl Let $X, Y, D, F$ be the same as in
Theorem 5.2. Let $f:\, X\rightarrow 2^Y\backslash\{\emptyset\}$ be a
set-valued map and $(x_0, y_0)\in {\rm gr} f$. Suppose that the
following conditions are satisfied:

{\rm (i)} \  for any $\preceq_F$-decreasing sequence $\{(x_n,
y_n)\}$ in $S_F(x_0, y_0)$, if $\{x_n\}$ is a Cauchy sequence, then
$\{x_n\}$ is convergent in $X$;

{\rm (ii)} \ $\xi$ {\rm (from (F3))} is lower bounded on
$P_Y(S_F(x_0, y_0))-y_0$ and $f(x)$ has the strict domination
property for any $x\in X$;

{\rm (iii)} \  for any $\preceq_F$-decreasing sequence $\{(x_n,
y_n)\}$ in $S_F(x_0, y_0)$, if $\{x_n\}$  converges to $x$, then
there exists $y\in f(x)$ such that $(x, y)\preceq_F (x_n, y_n),\
\forall n$.

Then, there exists $(\hat{x}, \hat{y})\in {\rm gr} f$ such that

{\rm (a) }  $y_0\in \hat{y} +F(\hat{x}, x_0) +D$ and $\hat{y}\in
{\rm SMin}^D f(\hat{x})$;

{\rm (b)}  for any $(x, y)\in {\rm gr}f \backslash \{(\hat{x},
\hat{y})\}$, $\hat{y}\not\in y+F(x, \hat{x}) +D.$}\\

{\bf Remark 5.2.}   \  Condition (i) in Corollary 5.1 can be
replaced by the following stronger condition

${\rm (i^{\prime})}$  for any sequence $\{(x_n, y_n)\}\subset {\rm
gr} f$, if $y_{n+1}\leq_D y_n, \ \forall n$, and $\{x_n\}$ is a
Cauchy sequence, then $\{x_n\}$ is convergent.

Compare (i$^{\prime}$) with the following $(f, D)$-lower
completeness (see Remark 3.1 or [41]): every Cauchy sequence
$\{x_n\}\subset X$ satisfying $f(x_n)\subset f(x_{n+1}) +D$ for
every $n$, is convergent. We see that condition (i$^{\prime}$) is
stronger than $(f, D)$-lower completeness. Let's call (i$^{\prime}$)
strong $(f,
D)$-lower completeness. Certainly, if $f$ is a vector-valued map, then the above two kinds of lower completeness are coincident.\\

For set-valued maps,  Khanh and Quy [31] introduced the following
concepts.\\

{\bf Definition 5.1.} \   Let $f:\, X\rightarrow
2^Y\backslash\{\emptyset\}$ be a set-valued map.

(i) $f$ is said to be $D$-lower semi-continuous from above(briefly,
denoted by $D$-lsca) at $\bar{x}$ if, for any convergent sequence
$x_n\rightarrow \bar{x}$ and any sequence $y_n\in f(x_n)$ with
$y_{n+1}\leq_D y_n\ \forall n$, there exists $\bar{y}\in f(\bar{x})$
such that $\bar{y}\leq_D y_n, \ \forall n.$

(ii)  $f$ is  said to be weak $D$-lower semi-continuous from above
(briefly, denoted by w.$D$-lsca) at $\bar{x}$ if, for each sequence
$x_n\rightarrow \bar{x}$ with $f(x_n)\subset f(x_{n+1}) +D.\ \forall
n$, one has $f(x_n)\subset f(\bar{x}) +D$.\\

As pointed out in [31], $D$-lsca implies w.$D$-lsca. We see that $f$
being w.$D$-lsca is exactly that $f$ is $D$-s.l.m.  Hence, $D$-lsca
maps can also be called strongly $D$-s.l.m. maps. If $f$ is a vector-valued map, then  $D$-s.l.m. and strong $D$-s.l.m. are coincident.\\

{\bf Corollary 5.2.} \ {\sl Let $(X, d)$ be a metric space, $Y$ be a
 locally convex space, $D\subset Y$ be a convex cone, $H\subset D$
be a convex set such that $0\not\in {\rm cl}(H+D)$ (i.e.,
$H^{+s}\cap D^+ \not= \emptyset$) and $H+D$ be $h_0$-closed for some
$h_0\in H$.

Let $f:\, X\rightarrow 2^Y\backslash \{\emptyset\}$ be strongly
$D$-s.l.m. (i.e., $D$-lsca) and $(x_0, y_0)\in {\rm gr} f$.

For $(x, y), (x^{\prime}, y^{\prime})\in {\rm gr} f$, define
$$(x^{\prime}, y^{\prime})\preceq_F( x, y)\ \ {\rm iff}\ \ y\in
y^{\prime} + d(x^{\prime}, x) H +D.$$ Suppose that the following
conditions are satisfied:

{\rm (i)} \    for any $\preceq_F$-decreasing sequence $\{(x_n,
y_n)\}$ in $S_F(x_0, y_0)$, if $\{x_n\}$ is a Cauchy sequence, then
$\{x_n\}$ is convergent in $X$ (particularly, $(X, d)$ is strongly
$(f, D)$-lower complete);

{\rm (ii)} \  there exists $\xi\in H^{+s}\cap D^+$ such that $\xi$
is lower bounded on $P_Y(S_F(x_0, y_0))$ and $f(x)$ has the strict
domination property for every $x\in X$.

Then the result of Corollary 5.1 holds.}\\

{\bf Proof.}\  It suffices to check assumption (iii) in Corollary
5.1. Let a sequence $\{(x_n, y_n)\}\subset {\rm gr}f$ satisfy
$(x_{n+1}, y_{n+1})\preceq_F (x_n, y_n)$ and $x_n \rightarrow
\bar{x}$. Clearly,  $y_{n+1}\leq_D y_n,\ \forall n$. By the
assumption, $f$ is strongly $D$-s.l.m., hence there exists
$\bar{y}\in f(\bar{x})$ such that $\bar{y}\leq_D y_n, \ \forall n$.
Next, we show that $(\bar{x}, \bar{y})\preceq_F (x_n, y_n), \
\forall n.$

If $d(x_n, \bar{x}) = 0$, then $x=x_n$.  In this case,
$\bar{y}\leq_D y_n$ is equivalent to that $y_n\in \bar{y} +
d(\bar{x}, x_n) H +D$. Certainly, we have $(\bar{x},
\bar{y})\preceq_F (x_n, y_n)$.

If $d(x_n, \bar{x}) >0$, take $i\in {\bf N}$ such that $d(x_n,
 \bar{x}) - (1/i)\,>\, 0$. Since $d(x_m, x_n) \rightarrow d(\bar{x},
 x_n)\ (m\rightarrow\infty)$, we may take $m>n$ such that $d(x_m,
 x_n) \geq d(\bar{x}, x_n)- (1/i)$. As $(x_m, y_m)\preceq_F (x_n,
 y_n)$ and $\bar{y}\leq_D y_m$, we have
\begin{eqnarray*}
y_n\,&\in&\, y_m + d(x_m, x_n) H +D\\
&\subset&\, y_m + (d(\bar{x}, x_n) -\frac{1}{i}) H +D\\
&\subset&\,\bar{y} + (d(\bar{x}, x_n) -\frac{1}{i})H+D.
\end{eqnarray*}
Thus,
\begin{eqnarray*}
y_n +\frac{1}{i} h_0\,&\in&\,\bar{y} + (d(\bar{x}, x_n)
-\frac{1}{i}) H +\frac{1}{i} H +D\\
&=&\, \bar{y} + d(\bar{x}, x_n) H +D.
\end{eqnarray*}
 From this,
 $$\frac{y_n-\bar{y} + ({1}/{i}) h_0}{d(\bar{x}, x_n)}\,\in\,
 H+D.$$
 Letting $i\rightarrow\infty$ and remarking that $H+D$ is
 $h_0$-closed, we have
 $$\frac{y_n-\bar{y}}{d(\bar{x}, x_n)}\,\in\, H+D\ \ {\rm and\ hence}\ \ y_n\in
 \bar{y} + d(\bar{x}, x_n) H +D.$$
 That is, $(\bar{x}, \bar{y})\preceq_F (x_n, y_n)$. \hfill\framebox[2mm]{}\\

As we have seen, the assumption that $\alpha H+D$ is closed for all
$\alpha >0$ (see [31, Theorem 5.2]) is not necessary. Here, we only
assume that $H+D$ is $h_0$-closed for some $h_0\in H$. In
particular, if $H$ is a singleton $\{k_0\}$ with $k_0\in
D\backslash-{\rm cl}(D)$, then we only need to assume that $D$ is
$k_0$-closed. As pointed out in [31, Section 5], the condition that
$f$ is strongly $D$-s.l.m. (i.e., $D$-lsca) and $f(x)$ has the
strict domination property imposed in Corollary 5.2 is equivalent to
the limiting monotonicity condition assumed in [5, Theorem 3.4].
Obviously, Corollary 5.2 improves [31, Theorem 5.2] and [5, Theorem
3.4]. Certainly, it also includes properly [4, Theorem 1].

As  done in the proof of [31, Theorem 5.3], from Corollary 5.2
we can obtain the following.\\

{\bf Corollary 5.3.} \ {\sl  Let $X,\, Y,\, D$ and $H$ be the same
as in Corollary 5.2 and additionally, $D$ be pointed and closed. Let
$f:\, X\rightarrow 2^Y\backslash\{\emptyset\}$ be $D$-lsc($D$-lower
semi-continuous), compact-valued and D-bounded (i.e., there exists a
bounded set $M$ such that $f(X)\subset M+D$), and $(x_0, y_0)\in
{\rm gr} f$. Suppose that assumption (i) of Corollary 5.2  is
satisfied (particularly, $(X, d)$ is strongly $(f, D)$-lower
complete). Then, the result of Corollary 5.2
holds.}\\

{\bf Proof.}\  Since $D$ is a closed convex pointed cone and $f$ is
compact-valued, $f$ has the strict domination property (see [31]).
As $f$ is D-bounded, for any $\xi\in H^{+s}\cap D^+$, $\xi$ is lower
bounded on $P_Y(S_F(x_0, y_0))$. Thus, assumption (ii) in Corollary
5.2 is satisfied. Also, $f$ being $D$-lsc and having compact-valued
implies that $f$ is $D$-lsca (see [31]). Now, we can
apply Corollary 5.2 and obtain the result. \hfill\framebox[2mm]{}\\

Concerning the strict domination property, there have been many
interesting results, for example, refer to [21, 43] and the
references therein. In fact, in Corollaries 5.1 and 5.2, the
condition that $f$ has the strict domination property can be
replaced by any one which implies $f$ having the strict domination
property. For example, from Corollary 5.2 we can obtain the
following Corollaries 5.4 and 5.5. First we recall some related
notions.

Let  $Y$ be a  locally convex space, $A\subset Y$  be nonempty,
$\Theta\subset Y$ be a bounded convex set and $D\subset Y$ be a
convex cone specifying a quasi-order $\leq_D$. Put
$\Theta_0:=\cup_{0\leq\lambda\leq 1} \lambda \Theta$. Then
$\cap_{\epsilon >0}(A-\epsilon \Theta_0)$ is called the
$\Theta$-closure of $A$ and denoted by ${\rm cl}_{\Theta}(A)$. If
${\rm cl}_{\Theta}(A)=A$, then $A$  is said to be $\Theta$-closed.
It is easy to see that $A$ is locally closed iff $A$ is
$\Theta$-closed for every  bounded convex set $\Theta$. And $A$ is
vectorial closed iff for every singleton $\Theta$, $A$ is
$\Theta$-closed. The following implications are obvious:
$${\rm closedness} \ \Longrightarrow\ {\rm local closedness}\
\Longrightarrow\ {\Theta}-{\rm closedness}.$$ But neither of two
converses is true, for details, see [43].

 A nonempty subset $A\subset Y$ is said to be $D$-complete (resp.,
$D$-locally complete) iff every Cauchy sequence (resp., locally
Cauchy sequence) $\{y_n\}\subset A$ with $y_{n+1}\leq_D y_n$, is
convergent (resp., locally convergent) to some point $\bar{y}\in A$.
It is easy to show that every $D$-complete set is $D$ locally complete.\\

{\bf Corollary 5.4.} \ {\sl Let $X,\, Y,\, D$ and $H$ be the same as
in Corollary 5.2 and, additionally, $D$ have a $\sigma$-convex base
$\Theta$.  Let $f:\, X\rightarrow 2^Y\backslash\{\emptyset\}$ be
strong $D$-s.l.m. (i.e., $D$-lsca), $D$-bounded and $f(x)$ be
$\Theta$-closed (particularly, locally closed or closed) for all
$x\in X$ and  $(x_0, y_0)\in {\rm gr} f$. Suppose that assumption
(i) of Corollary 5.2 is
satisfied. Then, the result of Corollary 5.2 holds.}\\

{\bf Proof.}\  Since $D$ has a $\sigma$-convex base $\Theta$ and
$f(x)$ is $\Theta$-closed, by [43, Corollary 5.2], $f(x)$ has the
strict domination property. Also, $f(X)$ being $D$-bounded implies
that assumption (ii)  in Corollary 5.2 is satisfied. Thus, the
result follows from Corollary 5.2. \hfill\framebox[2mm]{}\\

{\bf Corollary 5.5.} \ {\sl  Let $X,\, Y,\, D$ and $H$ be the same
as in Corollary 5.2 and, additionally, $D$ be locally closed and
have a bounded base. Let $f:\, X\rightarrow
2^Y\backslash\{\emptyset\}$ be strong $D$-s.l.m. (i.e., $D$-lsca),
$D$-bounded and $f(x)$ be $D$-locally complete (particularly,
$D$-complete or complete) for all $x\in X$ and $(x_0, y_0)\in {\rm
gr} f$. Suppose that assumption (i) of Corollary 5.2 is satisfied.
Then, the result of
Corollary 5.2 holds.}\\

{\bf Proof.}\  As done in the proof of Corollary 5.4, we only show
that $f(x)$ has the strict domination property. For any $x\in X$,
$f(x)$ is $D$-bounded and $f(x)$ is $D$-locally complete. Also, $D$
has a bounded base and $D$ is locally closed. By [43, Theorem 5.2],
$f(x)$ has the strict domination property and the proof is
completed.  \hfill\framebox[2mm]{}\\

Finally, as an application of Corollary 5.1, we give a Pareto
minimizer's version of Corollary 3.6.\\

{\bf Corollary 5.6.} \ {\sl Let $(X, d)$ be a metric space, $Y$ be a
locally convex space, $D\subset Y$ be a convex cone, $k_0\in
D\backslash -{\rm vcl}(D)$ and $D$ be $k_0$-closed. Let $f:\,
X\rightarrow 2^Y\backslash\{\emptyset\}$ be strongly s.l.m. (i.e.,
$D$-lsca) and $f(x)$ have the strict domination property for all
$x\in X$ and let $(X, d)$ be strongly $(f, D)$-lower complete.

Suppose the $(x_0, y_0)\in {\rm gr}f$ and $\epsilon >0$ such that
$y_0\not\in f(X) +\epsilon k_0 +D$.

Then for any $\lambda >0$, there exists $(\hat{x}, \hat{y}) \in {\rm
gr} f$ such that

{\rm (a)}  $y_0\in \hat{y} +(\epsilon/\lambda) d(\hat{x}, x_0) k_0
+D $ and $\hat{y}\in {\rm SMin}^D f(\hat{x})$;

{\rm (b)} for any $(x, y)\in {\rm gr} f\backslash\{(\hat{x},
\hat{y})\}$, $\hat{y}\not\in y+(\epsilon/\lambda) d(x, \hat{x}) k_0
+D$;

{\rm (c)} $d(x_0, \hat{x})\leq \lambda$.}\\

{\bf Proof.}\  We shall apply Corollary 5.1 to prove the conclusion.
Put
$$ F(x, x^{\prime}) := \,(\epsilon/\lambda) d(x, x^{\prime}) k_0,\ \
\forall x,\, x^{\prime}\in X.$$ Obviously, $F$ satisfies (F1) and
(F2). Since $k_0\in D\backslash -{\rm vcl}(D)$, we know that
$\xi_{k_0}(y)\not= -\infty,\ \forall y\in Y$. By Lemma 3.2,
$\xi_{k_0}$ is $D$-monotone and satisfies that
$$\xi_{k_0}(y+z) = \xi_{k_0}(y) + \xi_{k_0}(z),\ \ \forall y\in Y,\
\forall z\in F(X\times X).$$ Besides, for any $\delta >0$,
$$\zeta(\delta):=\,\inf\{\xi_{k_0}(y):\, y\in \cup_{d(x,
x^{\prime})\geq\delta}F(x, x^{\prime})\}\geq
(\epsilon/\lambda)\delta >0.$$ Hence, $F$ satisfies (F3) for
$\xi_{k_0}$. For $(x, y), \, (x^{\prime}, y^{\prime})\in {\rm gr}
f$, we define
$$(x^{\prime}, y^{\prime})\preceq_F (x, y)\ \Longleftrightarrow\
y\in y^{\prime} + (\epsilon/\lambda) d(x^{\prime}, x) k_0 +D.$$
Denote the set $\{(x, y)\in {\rm gr} f:\, (x, y) \preceq_F (x_0,
y_0)\}$ by $S_F(x_0, y_0)$.

Let $\{(x_n, y_n)\}$ be a $\preceq_F$-decreasing sequence in
$S_F(x_0, y_0)$ and let $\{x_n\}$ be a Cauchy sequence. From
$(x_{n+1}, y_{n+1})\preceq_F (x_n, y_n)$, we have $$y_n\in y_{n+1}
+(\epsilon/\lambda) d(x_{n+1}, x_n) k_0 +D \subset y_{n+1} +D\ \
{\rm and}\ \ y_{n+1}\leq_D y_n.$$ By the assumption that $(X, d)$ is
strongly $(f, D)$-lower complete, there exists $\bar{x}\in X$ such
that $x_n\rightarrow\bar{x} \ (n\rightarrow\infty)$, that is,
condition (i) in Corollary 5.1 is satisfied.

Let $\{(x_n, y_n)\}$ be a $\preceq_F$-decreasing sequence in
$S_F(x_0, y_0)$ and let $x_n\rightarrow \bar{x}$.  Since
$y_{n+1}\leq_D y_n$, $x_n\rightarrow \bar{x}$ and $f$ is strongly
$D$-l.s.m. (i.e., $D$-lsca), there exists $\bar{y}\in f(\bar{x})$
such that $\bar{y}\leq_D y_n,\ \forall n.$ For any given $n$, when
$m>n$, we have
\begin{eqnarray*}
y_n\,&\in&\, y_m + (\epsilon/\lambda) d(x_m, x_n) k_0 +D\\
&\subset&\, \bar{y} + (\epsilon/\lambda) d(x_m, x_n) k_0
+D.\hspace{7.9cm} (5.7)
\end{eqnarray*}
Since $d(x_m, x_n) \rightarrow d(\bar{x}, x_n)\
(m\rightarrow\infty)$ and $D$ is $k_0$-closed, from (5.7) we can
deduce that $y_n\in \bar{y} +(\epsilon/\lambda) d(\bar{x}, x_n) k_0
+D,$ and hence $(\bar{x}, \bar{y})\preceq_F (x_n, y_n)$. Thus,
condition (iii) in Corollary 5.1 is satisfied.

Since $y_0\not\in f(X) +\epsilon k_0 +D$, by Lemma 3.2, $\xi_{k_0}$
is lower bounded on $P_Y(S_F(x_0, y_0))-y_0$. And by the assumption
that $f(x)$ has the strict domination property, we know that
condition (ii) in Corollary 5.1 holds.

Now, applying Corollary 5.1, we obtain $(\hat{x}, \hat{y})\in {\rm
gr} f$ such that (a) and (b) hold. Finally, from (a), we have
$y_0\in \hat{y} + (\epsilon/\lambda) d(\hat{x}, x_0) k_0 +D$.
Combining this with $y_0\not\in f(X) +\epsilon k_0 +D$, we conclude
that $d(\hat{x}, x_0)\leq\lambda$. \hfill\framebox[2mm]{}\\

\noindent {\bf Acknowledgment}\\

The author is grateful to the referee and the editor for valuable
comments and suggestions which improved this paper.\\

\noindent{\bf References} \vskip 10pt
\begin{description}
\small
\item{[1]}  M. Adan, V. Novo, Proper efficiency in vector optimization
on real linear spaces, J. Optim. Theory Appl., 121 (2004), 515-540.

\item{[2]}  Y. Araya, Ekeland's variational principle and its
equivalent theorems in vector optimization, J. Math. Anal. Appl.,
 346 (2008), 9-16.

\item{[3]} J. -P. Aubin,  H. Frankowska, Set-Valued Analysis,
Birkh\"{a}user, Boston,  1990.

\item{[4]} T. Q. Bao, B. S. Mordukhovich, Variational principles for
set-valued mappings with applications to multiobjective
optimization, Control Cybern., 36 (2007), 531-562.

\item{[5]} T. Q. Bao, B. S. Mordukhovich, Relative Pareto minimizers
for multiobjective problems: existence and optimality conditions,
Math. Program, Ser.A, 122 (2010), 301-347.

\item{[6]}  E. M. Bednarczuk,  M. J. Przybyla, The vector-valued
variational principle in Banach spaces ordered by cones with
nonempty interiors, SIAM J. Optim. 18 (2007), 907-913.

\item{[7]}  E. M. Bednarczk, D. Zagrodny, Vector variational
principle, Arch. Math. (Basel), 93 (2009), 577-586.

\item{[8]} G. Y. Chen, X. X. Huang, A unified approach to the existing three
types of variational principle for vector valued functions, Math.
Methods Oper. Res., 48 (1998), 349-357.

\item{[9]}  G. Y. Chen, X. X. Huang, S. H. Hou, General Ekeland's
variational principle for set-valued mappings, J. Optim. Theory
Appl.,  106 (2000), 151-164.

\item{[10]}  G. Y. Chen, X. X. Huang, X. G. Yang, Vector Optimization, Set-Valued
and Variational Analysis, Springer-Verlag, Berlin, 2005.

\item{[11]} S. Dancs, M. Hegedus, P. Medvegyev, A general ordering
and fixed point principle in  complete metric space, Acta Sci. Math.
(Szeged), 46 (1983), 381-388.

\item{[12]} D. Dentcheva, S. Helbig, On variational principles,
level sets, well-posedness, and $\epsilon$-solutions in vector
optimization, J. Optim. Theory Appl., 89 (1996), 325-349.

\item{[13]}  I. Ekeland, Sur les prob\`{e}mes variationnels, C. R. Acad. Sci.
Paris 275 (1972), 1057-1059.

\item{[14]}  I. Ekeland,  On the variational principle, J. Math.
Anal. Appl. 47 (1974), 324-353.

\item{[15]} I. Ekeland, Nonconvex minimization problems, Bull. Amer.
Math. Soc. (N.S.) 1 (1979), 443-474.

\item{[16]}  C. Finet, L. Quarta, C. Troestler,  Vector-valued
variational principles, Nonlinear Anal. 52 (2003), 197-218.

\item{[17]}  F. Flores-Baz\'{a}n, C. Guti\'{e}rrez,  V. Novo,  A
Br\'{e}zis-Browder principle on partially ordered spaces and related
ordering theorems, J. Math. Anal. Appl., 375 (2011), 245-260.

\item{[18]} Chr. Gerstewitz (Tammer), Nichtkonvexe Dualit\"{a}t in
der Vektoroptimierung, Wiss. Z. TH Leuna-Merseburg 25 (1983),
357-364.

\item{[19]} Chr. Gerstewitz (Tammer), E. Iwanow, Dualit\"{a}t
f\"{u}r nichtkonvexe Vektoroptimierungsprobleme, Wiss. Z. TH llmenau
31  (1985), 61-81.

\item{[20]} Chr. Gerth (Tammer), P. Weidner, Nonconvex separation theorems and
some applications in vector optimization, J. Optim. Theory Appl., 67
(1990), 297-320.

\item {[21]}  A. G\"{o}pfert, H. Riahi,  Chr. Tammer,  C.
Z$\breve{a}$linescu,  Variational Methods in Partially Ordered
Spaces, Springer-Verlag, New York, 2003.

 \item{[22]} A. G\"{o}pfert,  C. Tammer and C. Z$\breve{a}$linescu,  On the
vectorial Ekeland's variational principle and minimal point theorems
in product spaces, Nonlinear Anal. 39 (2000),  909-922.

\item{[23]} C. Guti\'{e}rrez, B. Jim\'{e}nez, V. Novo, A set-valued
Ekeland's variational principle in vector  optimization, SIAM J.
Control. Optim., 47 (2008), 883-903.

\item{[24]} T. X. D. Ha, Some variants of the Ekeland  variational
principle for a set-valued map, J. Optim. Theory Appl., 124 (2005),
187-206.

\item{[25]}  A. H. Hamel, Equivalents to Ekeland's variational
principle in uniform spaces, Nonlinear Anal. 62 (2005), 913-924.

\item{[26]}  J. Horv\'{a}th, Topological Vector Spaces and Distributions, vol. 1,
Addison-Wesley, Reading, MA, 1966.

\item{[27]} G. Isac, The Ekeland's principle and the
Pareto $\epsilon$-efficiency, In: M. Tamiz (ed.) Multi-Objective
Programming and Goal Programming: Theories and Applications, Lecture
Notes in Econom. and Math. Systems, vol. 432, Springer-Verlag,
Berlin, 1996, 148-163.

\item{[28]} G, Isac, Chr. Tammer, Nuclear and full nuclear cones in
product spaces: Pareto efficiency and an Ekeland type variational
principle, Positivity, 9 (2005), 511-539.

\item{[29]}   P. Q. Khanh, D. N. Quy, On generalized Ekeland's variational
principle and equivalent formulations for set-valued mappings, J.
Glob. Optim., 49 (2011), 381-396.

\item {[30]} P. Q. Khanh, D. N. Quy, On Ekeland's variational
principle for Pareto minima of set-valued mappings, J. Optim. Theory
Appl.,  153 (2012), 280-297.

\item{[31]}  P. Q. Khanh, D. N. Quy, Versions of Ekeland's
variational principle involving set perturbations, J. Glob. Optim.
DOI 10.1007/s10898-012-9983-3.

\item{[32]}  D. Kuroiwa, On set-valued optimization, Nonlinear
Anal., 47 (2001), 1395-1400.

\item{[33]}  C. G. Liu, K. F. Ng, Ekeland's variational principle for
set-valued functions, SIAM J. Optim., 21 (2011), 41-56.

\item{[34]}  P. Loridan, $\epsilon$-Solutions in vector minimization problems, J.
Optim. Theory Appl., 43 (1984), 265-276.

\item{[35]}  A. B. N\'{e}meth, A nonconvex vector minimization problem,
Nonlinear Anal. 10 (1986),  669-678.

\item{[36]}  D. Pallaschke, S. Rolewicz, Foundations of Mathematical
Optimization, Math. Appl. 388, Kluwer, Dordrecht, 1997.

\item{[37]}  P. P\'{e}rez Carreras, J. Bonet, Barrelled Locally Convex Spaces,
North-Holland, Amsterdam, 1987.

\item{[38]} J. H. Qiu, Local completeness and dual local quasi-completeness,
Proc. Amer. Math. Soc. 129 (2001) 1419-1425.

\item{[39]} J. H. Qiu, Local completeness and drop theorem, J. Math.
Anal. Appl. 266 (2002), 288-297.

\item{[40]}  J. H. Qiu, A generalized Ekeland vector variational
principle and its applications in optimization, Nonlinear Anal., 71
(2009), 4705-4717.

\item{[41]} J. H. Qiu, On Ha's version of set-valued Ekeland's
variational principle, Acta Math. Sinica, English Series, 28 (2012),
717-726.

\item{[42]} J. H. Qiu, Set-valued quasi-metrics and a general
Ekeland's variational principle in vector optimization, SIAM J.
Control  Optim., 51 (2013), 1350-1371.

\item{[43]} J. H. Qiu, The domination property for efficiency and
Bishop-Phelps theorem in locally convex spaces, J. Math. Anal.
Appl., 402 (2013), 133-146.

\item{[44]} J. H. Qiu, F. He, A  general vectorial Ekeland's
variational principle with a p-distance, Acta Math. Sinica, 29
(2013), 1655-1678.

\item{[45]} J. H. Qiu, B. Li, F. He, Vectorial Ekeland's variational
principle with a w-distance and its equivalent theorems, Acta Math.
Sci., 32B (2012), 2221-2236.

\item{[46]} S. A. Saxon, L. M. S\'{a}nchez Ruiz, Dual local
completeness, Proc. Amer. Math. Soc., 125 (1997), 1063-1070.

\item{[47]} C. Tammer, A generalization of Ekeland's variational principle,
Optimization,  25 (1992)  129-141.

\item{[48]} C. Tammer, C. Z$\breve{a}$linescu, Vector variational
principle for set-valued functions, Optimization, 60 (2011),
839-857.

\item{[49]} M. Turinici, Maximal elements in a class of order
complete metric spaces, Math. Japon, 25 (1980), 511-517.

\item{[50]}  M. Turinici, Maximality principles and mean-value
theorems, Anais Acad. Brasil. Ciencias, 53 (1981), 653-655.

\item{[51]} C. Z$\breve{a}$linescu, Convex Analysis  in   General
Vector Spaces, World Sci., Singapore, 2002.

\end{description}
\end{document}